\documentclass[smallextended,final]{svjour3}       
\smartqed  
\usepackage{graphicx}
\usepackage{latexsym}
\usepackage{amsmath, amssymb}
\usepackage[color]{showkeys}
\usepackage{hyperref}
\definecolor{refkey}{gray}{.75}
\usepackage{booktabs}
\usepackage{enumerate}
\usepackage{xcolor}
\definecolor{rev}{RGB}{0,0,0}

\def\p{\partial}
\def\tilde{\widetilde}
\def\hat{\widehat}
\def\bs{\boldsymbol}

\def\ub{\boldsymbol{u}}
\def\vb{\bs{v}}

\def\taub{\bs{\tau}}
\def\kb{\mathbf{k}}
\def\nb{\mathbf{n}}
\def\tb{\mathbf{t}}

\def\CE{\mathcal{CE}}
\def\VE{\mathcal{VE}}
\def\VC{\mathcal{VC}}
\def\CV{\mathcal{CV}}
\def\EC{\mathcal{EC}}
\def\EV{\mathcal{EV}}
\def\VE{\mathcal{VE}}
\def\BC{\mathcal{BC}}
\def\IC{\mathcal{IC}}
\def\BE{\mathcal{BE}}
\def\IE{\mathcal{IE}}
\def\V{\mathcal{V}}

\def\M{\mathcal{M}}

\newcommand{\ra}[1]{\renewcommand{\arraystretch}{#1}}

\begin{document}

\title{Conservative finite volume schemes for the quasi-geostrophic equation on coastal-conforming  unstructured primal-dual meshes}

\titlerunning{Quasi-geostrophic}        

\author{Qingshan Chen         \and
        Lili Ju 
}

\authorrunning{Q.~Chen and L.~Ju} 

 \institute{Q.~Chen \at
               Department of Mathematical Sciences \\
               Clemson University\\
               \email{qsc@clemson.edu}           
            \and
            L.~Ju \at
            Department of Mathematics\\
            University of South Carolina\\
            \email{ju@math.sc.edu}           
 }

\date{\today}

\maketitle

\begin{abstract}
In this paper we propose finite volume schemes  for solving the inviscid and
viscous quasi-geostrophic equations on coastal-conforming unstructured
primal-dual 
meshes.  Several approaches for enforcing the boundary conditions are
also explored along with these schemes. The pure transport part in
these schemes are shown to conserve 
the potential vorticity along fluid paths in an averaged sense, and
conserve the potential enstrophy up to the time truncation
errors. Numerical tests  based on the centroidal Voronoi-Delaunay
meshes are performed to confirm these properties, and 
to distinguish the dynamical behaviors of these schemes. Finally some
potential  
applications of these schemes in different situations are discussed.

\keywords{Quasi-geostrophic equation; finite volume method; unstructured primal-dual meshes;
  vorticity-divergence formulation}
\end{abstract}

\section{Introduction}
\label{intro}

Due to the inherent instabilities of large-scale geophysical flows,
i.e.~the ocean and atmosphere, numerical simulations of these flows
cease to be a reliable reflection of the reality in the pointwise
sense after a few days. Then, it is natural that simulations of
large-scale geophysical flows set their goal on capturing the
long-term dynamics of these flows in the statistical sense. Numerical
schemes are thus developed with this goal in mind. The best known class of
schemes are the staggered-grid schemes of Arakawa and Lamb
(\cite{Arakawa1977-og}). Among this class of schemes, the most popular
and the most widely used one is the so-called C-grid scheme, for its
excellent representation of the inertial-gravity waves. The C-grid
scheme, also known as the MAC scheme in the classical CFD community
(\cite{Harlow1965-jc}), specifies the scalar quantities, such as the
mass and pressure, at cell centers, and the {\itshape normal} velocity
components at the cell edges. In the simulations, the tangential
velocity components are needed as well, due to the nonlinearity within
the models. In the original settings of rectangular meshes for the
C-grid or MAC schemes, 
the tangential velocity components can be easily reconstructed from
the normal velocity components. However, when the C-grid scheme is
implemented on unstructured meshes, the reconstruction becomes a
non-trivial matter, and is a significant source of errors (see
\cite{Ringler2010-sm,Thuburn2008-rc,Thuburn2014-ns}). One of the
authors of this present paper (Q.C.) made an attempt to remedy such
reconstruction issue by specifying both the normal and tangential
velocity components on the edges (\cite{Chen2013-fa}, and it has been
shown that, for the 
linearized shallow water equations with a constant Coriolis force, the
resulting scheme is decoupled on the primary and dual meshes
(\cite{Chen2016-eg}).  

Another type of numerical scheme that outperforms the C-grid in
resolving inertial-gravity waves is the non-staggered Z-grid scheme
(\cite{Randall1994-vu}). Indeed, it is clear that the Z-grid scheme
possesses the optimal dispersive wave relations among second-order
numerical schemes. The Z-grid scheme collocates the mass, vorticity,
and the divergence at the cell centers. Both the normal and tangential
velocity components can be recovered from the vorticity and divergence
by solving some Poisson equations
over the dual one of the primary-dual mesh. Thus, the tricky reconstruction
issue is avoided, at the cost of inverting some Laplace
operators. The Z-grid scheme has been implemented for the shallow
water equations on icosahedral meshes on a global sphere by Heikes and
Randall (\cite{Heikes1995-xx,Heikes1995-ky}). The cost issue
associated with inverting the Laplace operator on a global sphere was
dealt  by a multi-grid technique. However, we have
not seen any numerical and analytic studies of that scheme on bounded unstructured
meshes. 

We aim to adapt the Z-grid scheme for simulations of more general
oceanic flows over bounded domains. Our choice of the Z-grid scheme is
based on two factors. The first factor, an obvious one, is already
mentioned in above, that is, its optimal representation of the
dispersive wave relations. The second factor, which is equally
important, if not more so, is that the use of vorticity and divergence
as prognostic variables enable more accurate simulations of these
fundamental quantities of fluid dynamics. Due to the inherent
irregularities on unstructured meshes, the commonly used finite
difference and finite volume schemes on these meshes are only between
first and second-order accurate. For momentum-based numerical schemes,
this means that there is only a low, or even no accuracy for
quantities that are higher order derivatives of the momentum. 

In order to harness the great potentials of the Z-grid schemes,
several outstanding issues need to be addressed. The first issue
concerns the 
boundary conditions (BCs). For the momentum based systems and
schemes, the 
no-flux or no-slip BCs on the velocity are the natural choices. But
these natural 
choices do not have any counterparts on the vorticity or divergence variables. 
The second issue concerns the specifications of the diagnostic
variables to ensure the proper transport and/or conservation of certain
quantities such as the potential vorticity. Heikes and Randall
(\cite{Heikes1995-ky}) only partially solved this issue for the
icosahedral meshes and special types of flows. 
Lastly, in order for the vorticity-divergence based schemes to be
appealing for practical applications, its efficiency needs to be
improved, and be brought closer to, or on par with that of the
momentum-based numerical schemes. 

A large-scale geophysical flow such as the ocean or atmosphere is an
extremely complex dynamical system that exhibits a wide spectrum of
spatial and temporal scales and highly nonlinear phenomena. Developing
a reliable and 
efficient vorticity-divergence based numerical scheme for this system
demands careful planning and systematic efforts. Our plan is to
utilize a hierarchy of models that have already been proposed for
large-scale geophysical flows, and develop, implement, study, and test
vorticity-divergence based schemes for these models, starting from
simple models, and progressing to more complex and more realistic
ones. In this work, we deal with the barotropic quasi-geostrophic (QG)
equation. More 
complex models, such as the shallow water equations,
\textcolor{rev}{will be considered} 
in the subsequent works. 

The barotropic QG equation captures the dynamics of large-scale
geophysical flows in a single scalar equation that describes the
evolution of the QG potential vorticity (PV). It is the model
developed and used by Charney for the very first numerical
weather prediction on an electronic computer (\cite{Charney1950-sx}). 
The QG equation is already in a vorticity based form, and the velocity
field for the QG, as the leading component of the full velocity field, is
divergence free. For numerical simulation of the QG on an arbitrary
bounded domain, the most pressing
question is how to handle the boundary conditions. The question will
be one of  the main focuses of the current work. We will explore the proper setups
for the inviscid and viscous models. We will also explore the
difference and connection between artificial BC's, and implicitly and
explicitly enforced BC's. The QG equation is relatively simple
compared with other models, and therefore the specification or
computation of the diagnostics (e.g.~the stream function, the velocity
field) are mostly straightforward. Under the numerical setting and
discretization techniques that we choose in this paper, we are able
to show the pure advection part of our schemes conserve the potential
enstrophy up to time truncation error, and conserve the PV along the
fluid paths in an averaged sense. The third issue, namely efficiency,
will be left for future works, partly due to the fact that the
QG is already in a vorticity-based form, and there is not any
momentum-based numerical scheme for us to compare to.

We also note that a key part of the ingredients for our success in handling the BC's and the accuracy and 
conservation of the potential enstrophy and the PV is use of the
coastal-conforming centroidal Voronoi-Delaunay mesh that we
choose as the primal-dual mesh for the numerical simulations. Centroidal Voronoi tessellation (CVT) or spherical CVT (SCVT)  with their dual Delaunay triangulations have
been adopted for oceanic and atmospheric modeling by the MPAS (``Model for Prediction Across Scales'') team
(\cite{Ringler2013-pj,Ringler2008-rt,Ringler2011-lk,Ju2011-bx,Skamarock2012-zf}). The
C-grid staggering techniques are used, and 
the scalar quantities, such as mass, are specified at the Voronoi cell
centers, 
because the Voronoi cells, which are mostly hexagonal, are more
accommodating to the wave propagations compared with the dual
triangular cells (\cite{Danilov2010-il,Gassmann2011-if}). Meshes for
bounded ocean basins are culled out of meshes on the global sphere.
Meshes generated in this fashion are almost guaranteed
to have non-smooth boundaries: every connected boundary edges form an
angle of approximately $120^\circ$, as they are the edges of  (mostly)
hexagonal cells. This type of meshes will present an even more serious
complication for the vorticity-divergence based numerical schemes,
because the cell-centered scalar quantities such as the vorticity,
divergence, and stream functions, are always away from the
boundary. For these reasons, we take a different bottom-up
approach. Meshes will be generated for the ocean basins directly, with
the boundary fitted with partial Voronoi cells (or dual triangular
cells). With this approach, we are able to obtain a smooth
fitting to the coastal line, and have nodal points, where certain
scalar quantities are defined, on the boundary directly. On this mesh,
a complete discrete vector field theory has been developed
(\cite{Chen2016-gl}), which aides the subsequent development and
analysis of the numerical schemes. 

The rest of the paper is organized as follows. In Section \ref{s:cont}
we recall the barotropic QG equation and its theoretically required
boundary conditions and constraints. In Section \ref{mesh} the
specifications of the mesh are provided, and the main results of the
discrete vector field theory on this type of meshes are recalled. In
Section \ref{sec:scheme}, some finite volume schemes, for both the
inviscid and viscous models, are presented. 
Analyses of the transport part of the schemes concerning its
conservation of the PV and the potential enstrophy are carried out in
Section \ref{sec:anal}. 
Numerical tests and comparisons of these
schemes to demonstrate their qualitative behaviors are given in
Section \ref{s:numerics}. Concluding remarks and  some potential 
applications of these schemes in different situations are finally given  in Section
\ref{s:conclu}.  



\section{The continuous model and the boundary conditions}\label{s:cont} 
The barotropic QG equation reads
\begin{subequations}
  \label{eq:1}
  \begin{align}
    &\dfrac{\p}{\p t} q + \ub\cdot\nabla q =
    \dfrac{1}{H}\nabla\times{\taub} - \alpha\zeta  + 
                                            \mu\Delta\zeta,\label{eq:1a}\\
    &q = \zeta + \beta y - \dfrac{f_0}{H}(\psi-b),\label{eq:1b}\\
    &\ub = \dfrac{g}{f_0}\nabla^\perp \psi,\label{eq:1c}\\
    &\zeta = \nabla\times\ub \equiv \dfrac{g}{f_0}\Delta\psi.\label{eq:1d}
  \end{align}
\end{subequations}
In the above, $q$ stands for the potential vorticity, $\zeta$ the
relative vorticity, and $\ub$ the horizontal velocity field, and
$\psi$ the stream function. 
The parameter $f_0$ represents the average Coriolis parameter, $\beta$ (the beta term)
the meridional changing rate of the Coriolis parameter, $g$ the
gravity, $b$ the bottom topography, and $H$ the average fluid
depth. The friction with the bottom of 
the domain leads to an zeroth order damping to the vorticity field,
as revealed by the asymptotic analysis that produces the QG model
(\cite{Pedlosky1987-gk}). The strength of the damping is parameterized
by $\alpha$ on the right-hand side of \eqref{eq:1a}. The lateral
friction leads to a second-order damping to the vorticity field, in
the form of the Laplacian, as expected. The strength of the lateral
damping is parameterized by $\mu$. The vector field $\taub$ represents
the only external forcing to the system, which is the wind in the
physical world. 

The potential vorticity $q$ is a scalar quantity that is
being transported by the velocity field $\ub$, and they are related
through the stream function $\psi$ via \eqref{eq:1b},
\eqref{eq:1c} and \eqref{eq:1d}. In this sense, the potential 
vorticity can be considered an active tracer
(\cite{Friedlander2011-wb}). Within the specification \eqref{eq:1b}
for the PV, the zero order term $\psi$ represents the impact on the
vorticity from the fluctuations of the top surface. The fluctuations
of the top surface is usually small, and thus this term is often
omitted from the QG equation (\cite{Majda2006-qx}). However, in the
interior of the geophysical flows, where the interior Rossby
deformation radius is on the same order as the flow length scale, the
interface \textcolor{rev}{deforms} much more easily, and thus is a major
contributor to the vorticity (\cite{Pedlosky1987-gk}). For this
reason, we include the top 
surface fluctuation in our study. 

On a bounded domain $\M$, boundary conditions are required to complete the
system \eqref{eq:1}. When the lateral friction is absent, we consider
the no-flux boundary condition, which is the most natural choice for a
bounded domain surrounded by non-penetrable walls,
\begin{equation}
  \label{eq:bc1}
  \ub\cdot \nb = 0,\qquad {\rm on}\quad\p\M. 
\end{equation}
where $\nb$ represents the outer normal unit vector on the boundary
$\p\M$. 
In terms of the stream function $\psi$, this conditions means that
\begin{equation}
  \label{eq:bc1a}
  \dfrac{\p\psi}{\p \tb} = 0\quad \textrm{or}\quad \psi =
  \textrm{constant}, \qquad {\rm on}\quad\p\M, 
\end{equation}
where $\tb$ represents the positively oriented unit tangential vector
along the boundary $\p\M$. 
When the lateral friction is present, additional constraints are
required. In the context of a closed domain, the most natural choice
is the non-slip boundary condition,
\begin{equation}
  \label{eq:bc2}
  \ub\cdot\tb = 0, \qquad{\rm on}\quad \p\M.
\end{equation}
In terms of the stream function $\psi$, this conditions means that
\begin{equation}
  \label{eq:bc2a}
  \dfrac{\p\psi}{\p \nb} = 0, \qquad{\rm on}\quad \p\M. 
\end{equation}
A key step, in both the theoretical analysis and the numerical
computation of the QG equation is to recover the stream function $\psi$
from the potential vorticity $q$. The procedure for recovering the
streamfunction varies, depending on whether the friction is present or
not. 

\subsection{Without lateral friction}
When the  lateral friction is absent, the
stream function $\psi$ satisfies the elliptic equation \eqref{eq:1b}
and the no-flux BC's \eqref{eq:bc1a}. Compared with the standard
elliptic BVP with Dirichlet BC's, this system still has one parameter
to be determined, namely, the constant value of $\psi$ along the
boundary. We note that boundary value can not be arbitrarily set, due
to the zeroth order term on the right-hand side of \eqref{eq:1b}. 
An additional
constraint is needed in order to determine the constant value of
$\psi$ along the boundary. A logical option is the mass conservation
constraint,
\begin{equation}
  \label{eq:bc1b}
  \int_\M \psi dx = 0.
\end{equation}

Hence, in the inviscid case, the stream function $\psi$ is the
solution to a non-standard elliptic boundary value problem (BVP),
\begin{subequations}
  \label{bvp}
  \begin{align}
  & \zeta + \beta y - \dfrac{f_0}{H}(\psi-b) = q, &&\hspace{-1cm}{\rm in}\quad\M,\label{bvp1}\\
   &\psi = l,   &&\hspace{-1cm}{\rm on}\quad\p\M, \label{bvp2}\\
   &\int_\M \psi dx = 0.  &&\label{bvp3}
  \end{align}
\end{subequations}

The solvability of this non-standard elliptic BVP and the uniqueness
of its solution have been addressed 
in detail in \cite{Chen2017-fh}. Here, we recall the basic
strategy, since it has bearing on the numerical computation as well.  
We let $\psi_1$ and $\psi_2$ be solutions of the
following elliptic BVPs, respectively,
\textcolor{rev}{
\begin{subequations}
  \label{th10}
  \begin{align}
  \dfrac{g}{f_0}\Delta \psi_1 - \dfrac{f_0}{H}\psi_1 &= q - \beta y
         - \dfrac{f_0}{H}b, & &\hspace{-1cm}{\rm in}\quad\M,\label{th10a}\\
    \psi_1 &= 0, & &\hspace{-1cm}{\rm on}\quad\p\M,\label{th10b}
  \end{align}
\end{subequations}}
and
\textcolor{rev}{
\begin{subequations}
  \label{th11}
  \begin{align}
    \dfrac{g}{f_0}\Delta \psi_2 - \dfrac{f_0}{H}\psi_2 &= 0, &&\hspace{-1cm}{\rm in}\quad \M,\label{th11a}\\
    \psi_2 &= 1, &&\hspace{-1cm} {\rm on}\quad\p\M.\label{th11b}
  \end{align}
\end{subequations}}
By the standard elliptic PDE theories, both BVPs \eqref{th10} and
\eqref{th11} are well-posed under proper assumptions on the forcing
on the right-hand side of \eqref{th10a} and on the domain $\M$.
The solution to the original BVP \eqref{bvp} can be expressed in terms
of $\psi_1$ and $\psi_2$,
\begin{equation}
  \label{th12}
  \psi = \psi_1 + l\psi_2.
\end{equation}
The unknown constant $l$ can be determined using the mass conservation
constraint
\eqref{bvp3}
\begin{equation*}
  \int_\M \psi dx = \int_\M \psi_1 dx + l\int_\M \psi_2 dx =
  0,
\end{equation*}
which leads to 
\begin{equation}\label{th13}
l = -\dfrac{\displaystyle\int_\M \psi_1 dx}{\displaystyle\int_\M \psi_2 dx}. 
\end{equation}
We would like to point out that the expression \eqref{th13} for $l$ is valid
because, as a consequence of the maximum principle, $\psi_2$ 
is positive in the interior of the domain, and the integral $\psi_2$
in the denominator of \eqref{th13} is strictly positive.

The well-posedness of a frictionless (i.e.,~no bottom drag or lateral
diffusion) version of the model \eqref{eq:1} with the BC \eqref{bvp2}
and the constraint \eqref{bvp3} has been established in
\cite{Chen2017-fh}. 

\subsection{With lateral friction}
When the  lateral friction is present, both the
no-flux and the no-slip boundary conditions \eqref{eq:bc1a} and
\eqref{eq:bc2a} apply, however these conditions are over-specified for a
second-order elliptic equation \eqref{eq:1b}. To obtain a well-posed
system, one can apply the Laplacian operator $\Delta$ to \eqref{eq:1b}
to get a fourth-order biharmonic equation,
\textcolor{rev}{
\begin{equation}
  \label{eq:biharm}
     \dfrac{g}{f_0}\Delta^2\psi  - \dfrac{f_0}{H}\Delta \psi = q -
     \dfrac{f_0}{H}\Delta b.
\end{equation}}
The solution to the system \eqref{eq:biharm}, \eqref{eq:bc1a} and
\eqref{eq:bc2a} is not unique, since if $\psi$ is a solution, then so
is $\psi+c$ for any constant $c$. To fix this arbitrariness, one can
set $\psi=0$ on the boundary. Doing so will not affect the shape of
the stream function, and hence the velocity field, but it contradicts
the assumption of a free top surface. A better option is 
again to enforce the mass conservation constraint \eqref{eq:bc1b}. 

Hence, in the case with lateral \textcolor{rev}{friction}, one can determine
the stream function from the PV by solving a fourth-order BVP,
\begin{equation}
  \label{eq:bi-sys}\left\{
    \begin{aligned}
  &   \dfrac{g}{f_0}\Delta^2\psi  - \dfrac{f_0}{H}\Delta \psi = q -
     \dfrac{f_0}{H}\Delta b, & &{\rm in}\quad\M, \\
  & \dfrac{\p\psi}{\p \tb} = 0,& &{\rm on}\quad\p\M,\\
  & \dfrac{\p\psi}{\p \nb} = 0, & &{\rm on}\quad\p\M,\\
 & \int_\M \psi dx = 0. & &
    \end{aligned}\right.
\end{equation}

However, dynamically, this approach for determining the stream function
$\psi$ is defective. The beta term is a dominant term in the QG
dynamics. But, as a linear function, the beta term falls into the null
space of the Laplace operator $\Delta$. As a consequence, when the
Laplacian operator is applied to the original elliptic equation, the
beta term completely \textcolor{rev}{disappears}. In the next section, we will explore
alternative approaches for numerically determining the stream function
$\psi$ \textcolor{rev}{ when the lateral friction is present.}

\section{Mesh specification and the discrete vector field theory} 
\label{mesh}
\subsection{ Mesh specification}
\begin{figure}[!ht]
  \centering
  \includegraphics[width=4.5in]{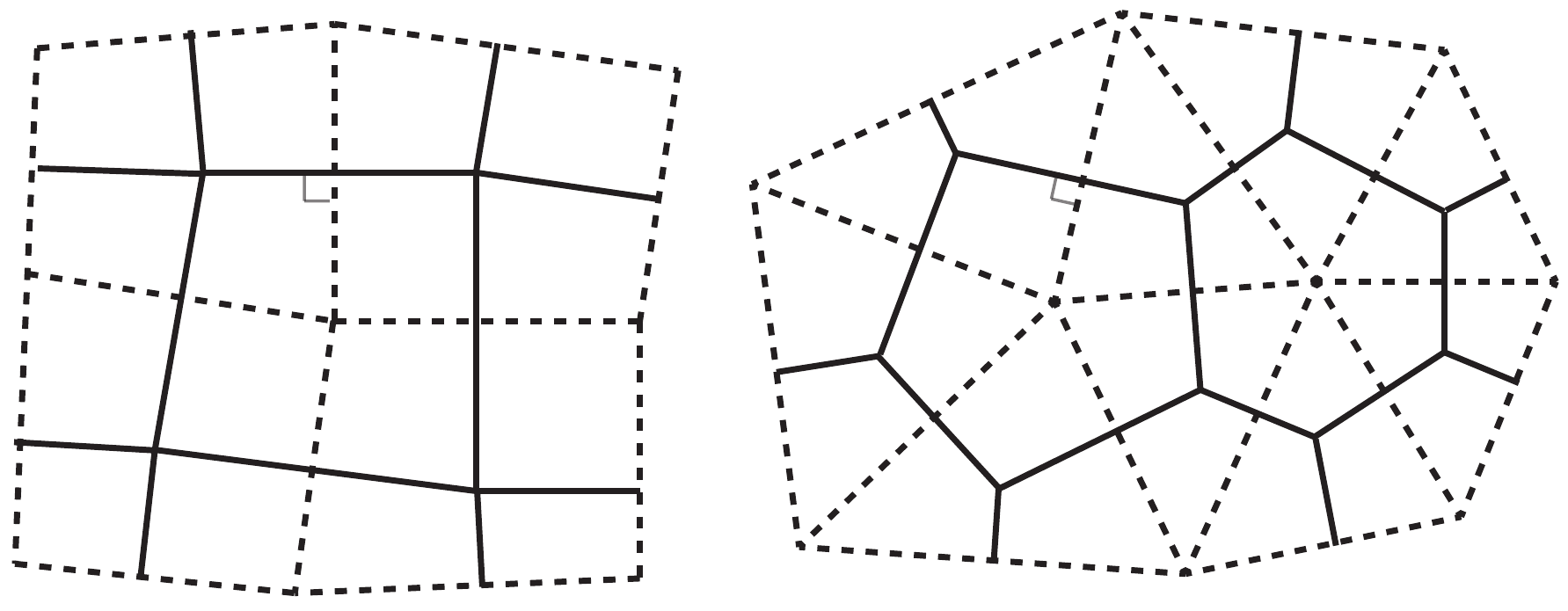}
  \caption{Generic conforming primal-dual meshes with the domain boundary passing
    through the primary cell centers, where the dark line denotes the primal mesh and the dashed line the dual mesh. Left: a generic quadrilateral
    dual mesh; right: a generic Voronoi-Delaunay mesh.}
  \label{fig:quad-dv}
\end{figure}

\textcolor{rev}{Our schemes for the QG equation  on based on very
  general orthogonal dual meshes, with one called {\em primal} and the
  other {\em dual}}.
The meshes consist of polygons, called cells, of arbitrary
shape. 
To avoid potential technical issues with the boundary, we 
assume that the domain $\M$ itself is also polygonal. 
The centers
of the cells on the primal mesh are the vertices of the cells on the
dual mesh, and vice versa. The edges of the primary cells intersect
{\it orthogonally} with the edges of the dual cells. The line segments
of the boundary $\partial\M$ pass through the centers of the
primary cells that border the boundary (Figure \ref{fig:p-edge}). 
\textcolor{rev}{Thus primary cells on the
boundary can be viewed as slices of regular cells like those in the
interior of the domain.}
{Two examples of this type of boundary conforming primal-dual meshes
  are shown in Figure \ref{fig:quad-dv}.}
  
  \begin{table}[!ht]
\centering
\ra{1.3}
\caption{Sets of basic mesh elements.}
\vspace{4mm}
\begin{tabular}{@{}l|l@{}}\toprule
Set & Definition\\
\midrule
$\IC$ & Set of interior cells\\
$\BC$ & Set of boundary cells\\
$\IE$ & Set of interior edges\\
$\BE$ & Set of boundary edges\\
$\V$ & Set of vertices\\
\bottomrule
\end{tabular}
\label{ta1}
\end{table}

\begin{figure}[!ht]
  \centering
  {\includegraphics[width=3.2in]{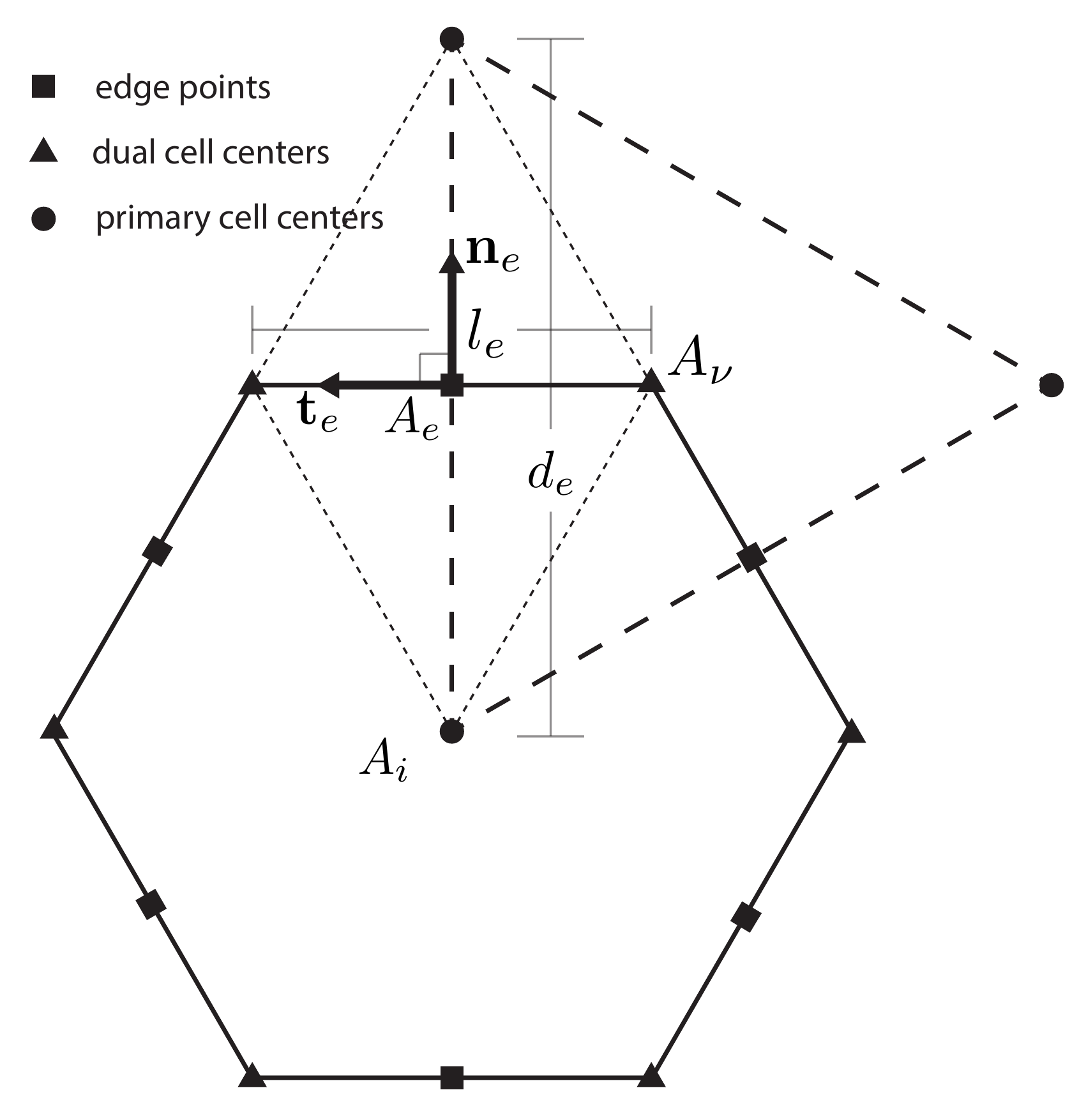}}
  \caption{Notations of a Voronoi-Delaunay mesh element.}
  \label{fig:notations}
\end{figure}

In order to effectively describe the schemes on these  meshes, some
notations are in order, for which we follow the conventions made in
\cite{Chen2013-fa,Ringler2010-sm}. The diagram in Figure
\ref{fig:notations} represents a  Voronoi-Delaunay mesh element. 
The primary cells are denoted as $A_i,\, 1\leq
i\leq N_c + N_{cb}$, where $N_c$ denotes the number of cells that are
in the interior of the domain, and $N_{cb}$ the number of cells that
are on the boundary. 
The dual cells, which
all lie inside the domain, are denoted as $A_\nu,\,1\leq \nu\leq
N_v$. 
The area of $A_i$ (resp.~$A_\nu$) is denoted as $|A_i|$
(resp.~$|A_\nu|$). 
Each primary
cell edge corresponds to a distinct dual cell edge, and vice
versa. Thus the primary and dual cell edges share a common index $e,\,
1\leq e\leq N_e+N_{eb}$, where $N_e$ denotes the number of edge pairs
that lie entirely in the interior of the domain, and $N_{eb}$ the
number of edge pairs on the boundary. Upon an edge pair $e$, the
distance between the two 
primary cell centers, which is also the length of the corresponding
dual cell edge, is denoted as $d_e$, while the distance between the
two dual cell centers, which is also the length of the corresponding
primary cell edge, is denoted as $l_e$. 
These two edges form the diagonals of a diamond-shaped region, whose
vertices consist of the two neighboring primary cell centers and the
two neighboring dual centers. The diamond-shaped region is also
indexed by $e$, and will be referred to as $A_e$.
The Euler formula for planar
graphs states that the number of primary cell centers $N_c + N_{cb}$, the
number of vertices (dual cell centers) $N_v$, and the number of
primary or dual cell edges $N_e + N_{eb}$ must satisfy the relation
\begin{equation}
\label{eq:g44}
  N_c + N_{cb} + N_v = N_e + N_{eb} + 1.
\end{equation}
The sets of basic elements (cells, edges, etc) are given in Table \ref{ta1}.
Additional  {sets of elements}, which  define the connectivity of the
primal-dual mesh, are given in Table \ref{ta2}. 

\begin{table}[!ht]
\centering
\ra{1.3}
\caption{Sets of elements defining the connectivity of a unstructured
  primal-dual mesh.}
\vspace{4mm}
\begin{tabular}{@{}l|l@{}}\toprule
Set & Definition\\
\midrule
$\EC(i)$ & Set of edges on the primary cell $A_i$\\
$\VC(i)$ & Set of dual cells on the primary cell $A_i$\\
$\CE(e)$ & Set of primary cells on the edge $e$\\
$\VE(e)$ & Set of dual cells on the edge $e$\\
$\CV(\nu)$ & Set of primary cells on the dual cell $D_\nu$\\
$\EV(\nu)$ & Set of edges on the dual cell $D_\nu$\\
\bottomrule
\end{tabular}
\label{ta2}
\end{table}

For each edge pair, a unit vector $\nb_e$, normal to the primary cell
edge, is specified. A second unit vector $\tb_e$ is defined as
\begin{equation}
\label{eq:g45}
  \tb_e = \kb\times\nb_e,
\end{equation}
with $\kb$ standing for the upward unit vector. Thus $\tb_e$ is
orthogonal to the dual cell edge, but tangent to the primary cell
edge, and points to the vertex on the left side of $\nb_e$. For each
edge $e$ and for each $i\in \CE(e)$ (the set of cells on edge $e$, see
Table \ref{ta2}), we define the direction indicator
\begin{equation}
\label{eq:g46}
  n_{e,i} = \left\{
  \begin{aligned}
    1& & &\phantom{sssssss}\textrm{if }\nb_e\textrm{ points away from primary cell
    }A_i,\\
    -1& &  &\phantom{sssssss}\textrm{if }\nb_e\textrm{ points towards primary cell
    }A_i,\\
  \end{aligned}\right.
\end{equation}
and for each $\nu\in \VE(e)$,
\begin{equation}
\label{eq:g47}
  t_{e,\nu} = \left\{
  \begin{aligned}
    1& & &\phantom{sssssss}\textrm{if }\tb_e\textrm{ points away from dual cell
    }A_\nu,\\
    -1& &  &\phantom{sssssss}\textrm{if }\tb_e\textrm{ points towards dual cell
    }A_\nu.\\
  \end{aligned}\right.
\end{equation}

For this study, we make the following regularity assumptions on the
meshes. We assume that the diamond-shaped region $A_e$ is actually
convex. In other words, the intersection point of each edge pair falls
inside each of the two edges. We also assume that the meshes are
quasi-uniform, in the sense that there exists $h>0$ such that, for 
each edge $e$,
\begin{equation}
\label{eq:g48}
  mh\leq l_e,\,d_e \leq Mh
\end{equation}
for some fixed constants $(m,\,M)$ that are independent of the
meshes. 
The staggered dual meshes are thus designated by $\mathcal{T}_h$.
{For the convergence analysis, it is assumed in
  \cite{Chen2016-gl} that, for each edge pair $e$, the primary cell
  edge nearly bisect 
  the dual cell edge, and miss by at most $O(h^2)$. This assumption is
  also made here for the error analysis. Generating meshes
  satisfying this requirement on irregular domains, i.e.~domains
  with non-smooth boundaries or domains on surfaces, can be a
  challenge, and will be addressed elsewhere. But we point out that, on
  regular domains with smooth boundaries, this type of meshes can be
  generated with little extra effort in addition to the use of
  standard mesh generators, such as the 
  centroidal Voronoi tessellation (CVT) or Spherical CVT (SCVT)  meshing algorithm (\cite{Du2003-gn,Du2002-lf,Du1999-th})}.  

\subsection{Specification of the discrete differential
  operators}\label{sec:discrete-ops} 
In this section, we denote by $\psi_h$ a discrete scalar field defined
at the primary cell centers, and by $\tilde\varphi_h$ a discrete  scalar field
defined at the primary cell vertices,
\begin{align*}
  \psi_h &= \sum_{i\in \IC\cup\BC} \psi_i \chi_i,\\
  \tilde\varphi_h &= \sum_{\nu\in \IC\cup\BC} \tilde\varphi_\nu \chi_\nu.
\end{align*}
Here, $\chi_i$ and $\chi_\nu$ stand for the piecewise constant
characteristic functions on the primary cells and on the dual cells,
respectively. 
A discrete scalar field on the primary mesh can be mapped to a
discrete scalar field on the dual mesh via area-weighted averaging,
\begin{equation}
  \label{eq:remap}
  \tilde\psi_\nu = \tilde{[\psi_h]}_\nu \equiv
  \dfrac{1}{|A_\nu|}\sum_{i\in \CV(\nu)} \psi_i |A_{i,\nu}|. 
\end{equation}
In the above, $|A_{i,\nu}|$ stands for the area of the kite-shaped
region between cell $i$ and vertex $\nu$ (see Figure
\ref{fig:notations}). 
A discrete scalar field on the primary mesh can also be mapped to a
discrete scalar field on the edges via arithmetic mean,
\begin{equation}
  \label{eq:c2e}
  \hat \psi_e = \hat{[\psi_h]}_e \equiv = \dfrac{1}{2}\sum_{i\in
    \CE{e}} \psi_i.
\end{equation}

\begin{figure}[h]
  \centering
  \includegraphics[width=3.5in]{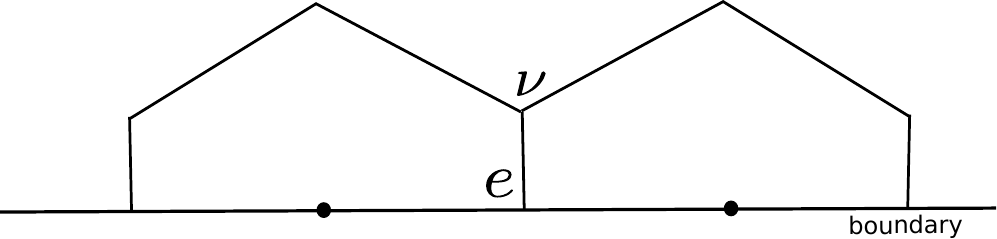}
  \caption{\textcolor{rev}{
With each (partial) edge on the boundary, only one vertex
    is associated.}}
  \label{fig:p-edge}
\end{figure}

On each edge $e$, the
discrete gradient operator on $\psi_h$
is defined as 
\begin{equation}
  [\nabla_h\psi_h]_e = \dfrac{-1}{d_e}\sum_{i\in \CE(e)}\psi_i
  n_{e,i},\label{eq:d34} 
 \end{equation}
and the skewed discrete gradient operator on $\tilde\varphi_h$ 
is defined as
\begin{equation}
  [\tilde\nabla_h^{\perp}\tilde\varphi_h]_e = \dfrac{1}{l_e}\sum_{\nu\in
    \VE(e)}\tilde\varphi_\nu t_{e,\nu}.\label{eq:d35} 
 \end{equation}
 With
each boundary edge $e$, since only one vertex is associated (Figure
\ref{fig:p-edge}), the 
definition \eqref{eq:d35} can be written as 
\begin{equation}
  [\tilde\nabla_h^{\perp}\tilde\varphi_h]_{e\textrm{ on boundary}} =
  \dfrac{1}{l_e}\tilde\varphi_\nu t_{e,\nu},\label{eq:d40} 
 \end{equation} 
where $\nu$ is the single element in $\VE(e)$. This amounts to
implicitly requiring that $\tilde\varphi_h$ vanishes on the boundary. 
We define the gradient of $\psi_h$ and the skew gradient of
$\tilde\varphi_h$ as
\begin{align} 
  &\nabla_h \psi_h = \sum_{e=1}^{N_e + N_{eb}} [\nabla_h\psi_h]_e \chi_e
  \nb_e,\label{eq:d36}\\ 
  &\tilde\nabla_h^{\perp} \tilde\varphi_h = \sum_{e=1}^{N_e + N_{eb}}
  [\tilde\nabla_h^{\perp}\tilde\varphi_h]_e \chi_e \nb_e.\label{eq:d37} 
\end{align}
Here $\chi_e$ represents the piecewise constant characteristic
\textcolor{rev}{function} on the diamond-shaped regions $A_e$. 

We denote by $\ub_h$ and $\vb_h$ discrete vector fields
defined at the edges, in the normal (w.r.t.~the primary cells) and the
tangential directions, respectively,
\begin{align*}
  \ub_h &= \sum_{e\in\IE\cup\BE} u_e \nb_e \chi_e,\\
  \vb_h &= \sum_{e\in\IE\cup\BE} v_e \tb_e \chi_e.
\end{align*}

Around each primary cell $i$, a discrete divergence operator can be
defined on the discrete vector field $\ub_h$,
\begin{equation}
  \label{eq:d26}
  \left[\nabla_h \cdot \ub_h\right]_i = \dfrac{1}{A_i}\sum_{e\in \EC(i)}u_e l_e n_{e,i}.
\end{equation}
It is worth noting that, on partial cells on the boundary, the
summation on the right-hand side only includes fluxes across the edges
that are inside the domain and the partial edges that intersect with
the boundary, and this amounts to imposing a no-flux
condition across the boundary, i.e. $\ub\cdot\nb = 0$ on $\p\M$.
It is clear from the definition \eqref{eq:d26} that the image of the discrete divergence operator
$\nabla_h\cdot$ on $\ub_h$ is a scalar field on the primary mesh, which can be written as
\begin{equation}
  \label{eq:d27}
  \nabla_h\cdot \ub_h = \sum_{i\in\IC\cup\BC} \left[\nabla_h \cdot
    \ub_h\right]_i\chi_i.
\end{equation}
Around each dual cell $\nu$, a discrete curl operator can be defined
on $\ub_h$,
\begin{equation}
\label{eq:d28}
  \left[\tilde\nabla_h \times \ub_h\right]_\nu = \dfrac{-1}{A_\nu}\sum_{e\in \EV(\nu)}u_e d_e t_{e,\nu}.
\end{equation}
The tilde atop $\nabla$ again signifies the involvement of the dual cells.
Thus, the image of the discrete curl operator
$\tilde\nabla_h\times$ on each $\ub_h$ is a scalar field defined on
the dual mesh,
\begin{equation}
\label{eq:d29}
  \tilde\nabla_h\times \ub_h = \sum_{\nu=1}^{N_v}\left[\tilde\nabla_h \times
    \ub_h\right]_\nu\chi_\nu.
\end{equation}

For a scalar field $\psi_h$ defined at cell centers, the discrete
Laplacian operator $\Delta_h$ can also be defined,
\textcolor{rev}{\begin{equation}
  \label{eq:d30}
  \Delta_h \psi_h := \nabla_h\cdot\left(\nabla_h\psi_h\right) 
\end{equation}}

\textcolor{rev}{
Finally, in this study, we need an inner product for scalar fields on
the primary mesh, and an inner product for discrete normal vector
fields on the edges,
\begin{align}
  (\psi_h,\,\varphi_h)_{0,h} =& \sum_{i\in\IC\cup\BC} \psi_i\varphi_i
                          A_i,\label{eq:inp1}\\
  (\ub^1_h,\,\ub^2_h)_{0,h} =& \sum_{e\in \IE\cup\BE} u^1_e u^2_e
                               A_e.\label{eq:inp2} 
\end{align}}

The discrete forms of the integration by parts formulas are also available.
\begin{lemma}\label{lem:integ-by-parts}
  For the discrete normal vector field $\ub_h$ and the discrete
  scalar fields $\psi_h$ and $\tilde\varphi_h$, the
  following relations hold:
  \begin{align}
    \label{eq:d41}
    \left(\ub_h,\,\nabla_h\psi_h\right)_{0,h} &=
    -\dfrac{1}{2}\left(\nabla_h\cdot \ub_h,\,\psi_h\right)_{0,h},\\
    \left(\ub_h,\,\tilde\nabla_h^\perp\tilde\varphi_h\right)_{0,h} &=
    -\dfrac{1}{2}\left(\tilde\nabla_h\times
    \ub_h,\,\tilde\varphi_h\right)_{0,h}.\label{eq:d42} 
  \end{align}
\end{lemma}
The proof can be found in \cite{Chen2016-gl}. \textcolor{rev}{We
  simply point out that 
these integration by parts formulas hold thanks to the no-flux
boundary condition on $\ub_h$, implied in the specification of the
discrete divergence operator $\nabla_h\cdot(\,)$, and the
homogeneous Dirichlet boundary condition on $\tilde\varphi_h$, implied
in the specification of the skewed gradient operator
$\tilde\nabla_h^\perp$. }

\section{Finite volume schemes for the QG equation}\label{sec:scheme}
\subsection{Finite volume schemes for the inviscid QG equation} 

We deal with the system \eqref{eq:1}, but with the viscosity
$\mu$ set to zero. The proper boundary conditions for the resulting
system are the no-flux boundary condition \eqref{eq:bc1a} and the mass
conservation constraint \eqref{eq:bc1b}. 

\subsubsection{A finite volume scheme for the inviscid QG equation (IVFV1)}
\label{s:fv-1}
The key quantity for the QG equation is the PV $q$. We approximate
this quantity with a discrete scalar field $q_h$ that is defined at
cell centers, i.e.,
\begin{equation}
  \label{eq:qh}
  q_h = \sum_{i\in \IC\cup \BC} q_i \chi_i. 
\end{equation}

Closely related to the PV is the discrete stream function $\psi_h$,
defined as 
\begin{equation}
  \label{eq:psih}
  \psi_h = \sum_{i\in\IC\cup\BC} \psi_i\chi_i. 
\end{equation}
One can directly derive the stream function from the PV by solving 
a discrete elliptic BVP, which is a FV discretization of the continuous
BVP \eqref{bvp},
\begin{equation}\label{eq:disc-e}
  \left\{
    \begin{aligned}
      &\dfrac{g}{f_0} [\Delta_h \psi_h]_i - \dfrac{f_0}{H} \psi_i = q_i
      - \beta y - \dfrac{f_0}{H}b_i, & &\quad i\in \IC,\\
      &\psi_i = \textrm{constant} & &\quad i\in \BC,\\
      &\sum_{i\in \IC\cup\BC} \psi_i A_i = 0. & &
    \end{aligned}\right.
\end{equation}
We note that the BC \eqref{bvp2} and the constraint \eqref{bvp3} are
enforced explicitly. \textcolor{rev}{This discrete system can be
  solved using the procedure outlined in \eqref{th10}--\eqref{th13}.}

The velocity field $\ub$ in the QG system is non-divergent, thanks 
to \eqref{eq:1c}. As a result, the nonlinear advection term in the
transport equation \eqref{eq:1a} can be written as a net flux of the
PV $q$. With the discrete PV $q_h$ defined at the cell centers, the
most natural way to compute its fluxes is to define the normal
velocity components at cell edges and use the discrete divergence
theorem. According to \eqref{eq:1c}, the normal velocity component $u_e$ can be
computed using the stream function defined at the cell
vertices. 
Therefore, to compute the velocity field $\ub_h$, we first map $\psi_h$
to a scalar field 
$\tilde\psi_h$ on the cell vertices via \eqref{eq:remap},
\begin{equation}
  \label{eq:psi-tilde}
  \tilde\psi_h = \tilde{[\psi_h]}. 
\end{equation}
 The discrete velocity field is then computed using the discrete
skewed gradient operator defined in \eqref{eq:d37}
\begin{subequations}
\label{eq:u_h}
\begin{align}
  \ub_h &= \sum_{e\in \IE\cup \BE} u_e \chi_e\nb_e  \label{eq:uh}\\
  &= \tilde\nabla_h^\perp \tilde\psi_h.   \label{eq:uh-def}
\end{align}
\end{subequations}

\begin{figure}[!ht]
  \centering
    \includegraphics[width=2.5in]{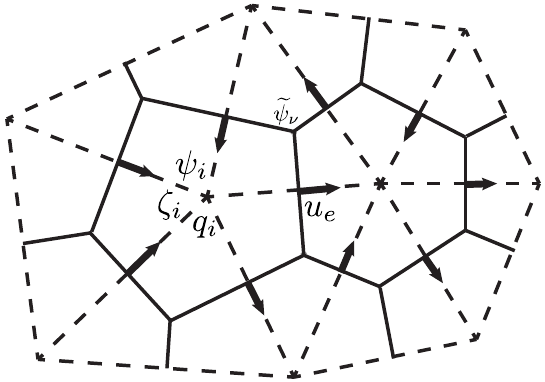}
  \caption{Placement of the discrete variables.}
  \label{fig:1}
\end{figure}

The placement of these discrete variables is shown on Figure
\ref{fig:1} in details. We recall that the PV $q_h$ is the only prognostic
variable, and the rest are all diagnostic variables, and can be
computed (or diagnosed) from the prognostic variables. 
The discrete PV filed evolves according to the following semi-discrete
equation,
\begin{equation}
\label{eq:q-evol}
  \dfrac{d q_i}{d t} + \left[ \nabla_h \cdot(\ub_h \hat q_h)\right]_i =
  \dfrac{1}{H}[\nabla\times\taub]_i - \alpha \zeta_i,\qquad
  i\in\IC\cup\BC. 
\end{equation}
\textcolor{rev}{In the above, $[\nabla\times\taub]_i$ represents the value of the curl
of the external forcing (wind) on cell $i$, which is calculated from
the given function or data,  and $\hat q_h$, the discrete PV field at the cell edges, is
a diagnostic variable computed from $q_h$ by the cell-to-edge mapping}
\eqref{eq:c2e}, 
\begin{equation}
\label{eq:qhat_e}
  \hat q_h = \hat{[q_h]} = \sum_{e\in\IE\cup\BE} \hat q_e \chi_e. 
\end{equation}

The discrete relative vorticity is collocated with the PV $q_h$ and
the stream function $\psi_h$. It can be computed from the PV $q_h$,
using the equation \eqref{bvp1},
\begin{subequations}\label{eq:zeta_h}
\begin{align}
  \zeta_h &= \sum_{i\in \IC\cup\BC} \zeta_i \chi_i,  \label{eq:zeta-i-1}\\
   \zeta_i &= q_i - \beta y_i + \dfrac{f_0}{H}\left(\psi_i -
  b_i\right).\label{eq:zeta_h-def}
\end{align}  
\end{subequations}

The advantage of this scheme is that the theoretical BC's \eqref{bvp2}
and \eqref{bvp3} are enforced explicitly, and no artificial BC's are
enforced or implied. The downside of this scheme is that, by
discretizing the elliptic equation \eqref{bvp1} up to the boundary, it
assumes that it holds up to the boundary, which is only possible if the
stream function is smooth over the entire domain
(i.e.~$C^2(\overline\M)$). The last requirement is unusually high in
theoretical analysis of this type of problems. 

In summary, the FV scheme IVFV1 for the inviscid QG equation
consists of the semi-discrete equation \eqref{eq:q-evol} and the
auxiliary equations \eqref{eq:disc-e}, \eqref{eq:psi-tilde},
\eqref{eq:u_h},  
\eqref{eq:qhat_e}, and \eqref{eq:zeta_h}. 

\subsubsection{A FV scheme with artificial boundary conditions for the
  inviscid QG equation (IVFV2)}\label{s:fv-2}
\begin{figure}[!ht]
  \centering
  \includegraphics[width=9cm]{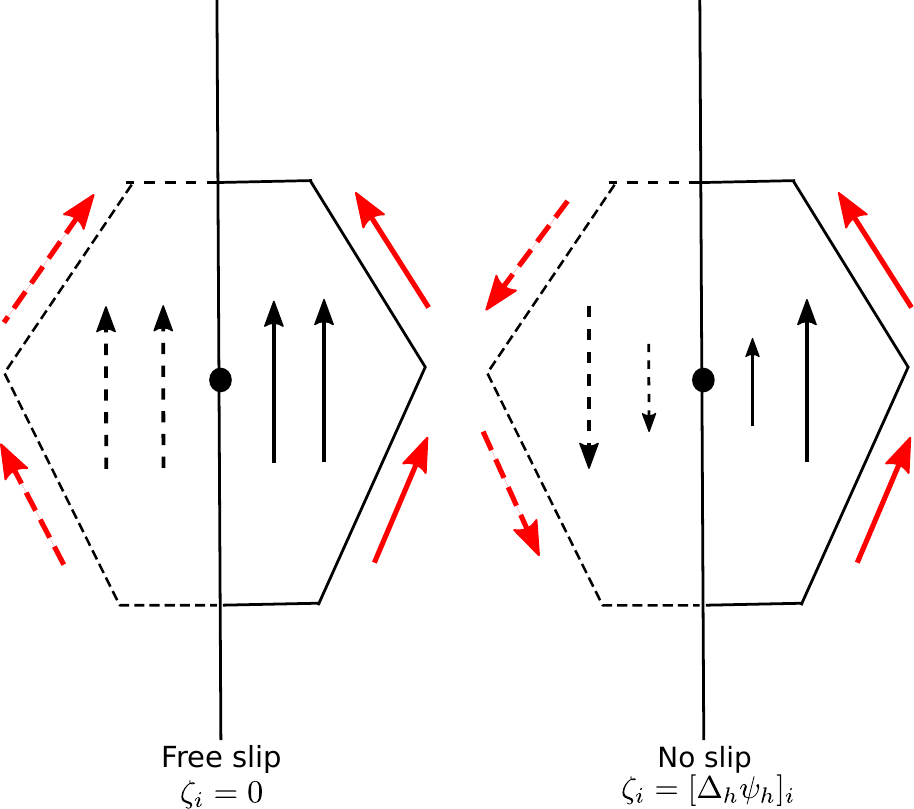}
  \caption{Boundary conditions on the vorticity.}
  \label{fig:vort-bcs}
\end{figure}

The previous numerical scheme IVFV1 leaves the PV $q_h$ freely
advected by the 
velocity field, even on the boundary. The boundary conditions are
applied on the stream function only, and they affect the evolution of
the PV $q_h$ through the stream function. No explicit constraints are
imposed on the boundary values of the PV. To maintain some control
on the boundary behavior of the PV, an artificial BC can be
utilized. For the inviscid flow and {\itshape in the absence of
  external forcing 
near the boundary}, the fluids there will flow along the boundary, and
there should be no velocity differential in the immediate vicinity of
the boundary. We consider a partial cell on the boundary (Figure
\ref{fig:vort-bcs}). We 
imagine that the partial cell is completed by a ghost partial cell
(dashed) through mirroring. Then the cell center of the original
partial cell becomes the center of the whole cell. Within this cell
and near the boundary,
due to the constraint exerted by the rigid but frictionless boundary, the
flow has a uniform velocity along the direction of the 
boundary. Thus, at the cell center, the vorticity should be zero
identically. The same conclusion will result when one project the
velocity field onto the edges, and then apply Green's theorem on the
tangential velocity components around the cell (red arrows). Thus, we
propose the following BC for the relative vorticity,
\begin{equation}
  \label{eq:vort-bc}
  \zeta_i = 0,\qquad i\in \BC.
\end{equation}
For the PV, this condition implies that
\begin{equation}
  \label{eq:pv-bc}
  q_i = \beta y_i - \dfrac{f_0}{H}\left(\psi_i - b_i\right),\qquad
  i\in \BC.
\end{equation}

Physically, the condition \eqref{eq:pv-bc} means that, on the
boundary, the PV equals the local planetary vorticity plus the
variations in the top surface and the bottom topography. Both the
planetary vorticity and the bottom topography are fixed in time; only
the top surface is time dependent.

\begin{figure}[!ht]
  \centering
  \includegraphics[width=3.5in]{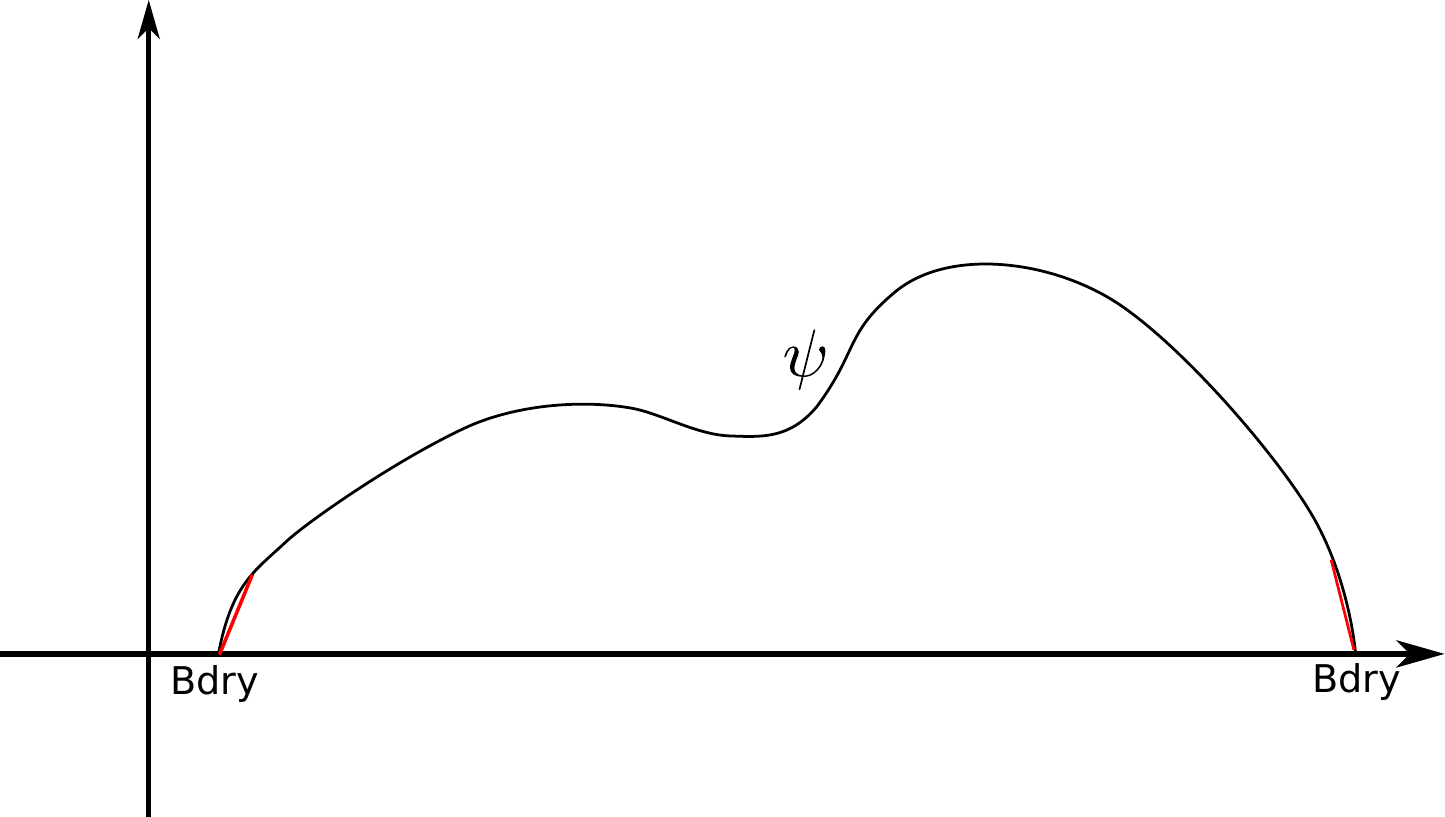}
  \caption{Artificial boundary conditions on $\psi$.}
  \label{fig:psi-bc}
\end{figure}

Near the boundary, the velocity component along the boundary is given
by $\p\psi/\p \nb$, and the vorticity is given by $\p^2\psi/\p \nb^2$,
where $\nb$ designates the outer normal unit vector on the
boundary. Thus, in terms of the continuous variables, the artificial
BC \eqref{eq:vort-bc} means that
\begin{equation}
  \label{eq:psi-bc-2}
  \dfrac{\p^2\psi}{\p \nb^2} = 0,\qquad {\rm on}\quad\p\M.
\end{equation}
Geometrically, this condition mandates that, near the boundary, the
stream function be linear along the direction normal to the boundary
(Figure \ref{fig:psi-bc}). \textcolor{rev}{\itshape We point out that
  this BC is unphysical, not needed by 
the continuous system, and it cannot be supported either.} But its
discrete version, namely \eqref{eq:vort-bc} and \eqref{eq:pv-bc}, can
help to exert a certain amount of control on the behavior of the PV on
the boundary.

The rest of the scheme will remain the same as the previous scheme,
except that the scope of the indices for certain equations
\eqref{eq:q-evol} need to be adjusted to reflect the newly imposed
BC's. Specifically, the discrete equation \eqref{eq:q-evol} is changed
to 
\begin{equation}
\label{eq:q-evol-in}
  \dfrac{d q_i}{d t} + \left[ \nabla_h \cdot(\ub_h \hat q_h)\right]_i =
  \dfrac{1}{H}[\nabla\times\taub]_i - \alpha \zeta_i,\qquad
  i\in\IC,
\end{equation}
and the equation \eqref{eq:zeta_h-def} is changed to
\begin{equation}
  \label{eq:zeta_i-in}
    \zeta_i = q_i - \beta y_i + \dfrac{f_0}{H}\left(\psi_i -
  b_i\right),\qquad i\in\IC
\end{equation}
for computing the vorticity in interior cells. 

In summary, the FV scheme with an artificial boundary condition for the
inviscid QG equation
consists of the semi-discrete equation \eqref{eq:q-evol-in} and the
auxiliary equations \eqref{eq:disc-e}, \eqref{eq:psi-tilde},
\eqref{eq:u_h},  \eqref{eq:vort-bc}, \eqref{eq:pv-bc},
\eqref{eq:qhat_e}, and \eqref{eq:zeta_i-in}.

\subsection{Finite volume schemes for the viscous QG equation} 

The viscous QG equation \eqref{eq:1} with a non-zero viscosity $\mu$ requires
two boundary conditions \eqref{eq:bc1a} and \eqref{eq:bc2a} plus one
constraint $\eqref{eq:bi-sys}_4$. It has been pointed out in Section
\ref{s:cont} that, if one wants to determine the stream function $\psi$
from the PV $q$, one needs to differentiate the elliptic equation
\eqref{eq:1b} to obtain a fourth-order biharmonic equation so that the
full boundary conditions can be enforced. However, doing so will kill
the beta term in the specification of the PV $q$, and will have a
detrimental impact on the dynamical behavior of the system. Numerical
tests have demonstrated that this approach leads to an
over-diffusive system, and that the phenomenon that relies on the beta
effect, namely the intensification of the western boundary current, is
largely gone. To remedy this issue, we propose two different approaches.

\subsubsection{A semi-implicit finite volume scheme for the viscous QG
  equation (VSFV1)}\label{s:fv-3}
 The first approach is to skip the intermediate
step on the PV advection, and solve for the stream function from a
semi-implicit discretization of the fourth-order equation
\eqref{eq:1a} directly. As before, the stream function $\psi$ is
approximated by a discrete scalar field $\psi_h$ defined at the cell
centers. From the discrete stream function $\psi_h$, the discrete relative
vorticity $\zeta_h$ can be full expressed on the interior cells, 
\begin{equation}
  \label{eq:zeta-ic}
  \zeta_i = [\Delta_h \psi_h]_i,\qquad i\in\IC.
\end{equation}
The PV $q_h$ is also fully expressed on the interior cells, according
to \eqref{eq:1b} and \eqref{eq:1d},
\begin{equation}
  \label{eq:q-ic}
  q_i = [\Delta_h\psi_h ]_i + \beta y_i - \dfrac{f_0}{H} (\psi_i -
  b_i), \qquad i\in\IC.
\end{equation}
The diffusion $\Delta\zeta$ can be expressed on the interior cells
that are not connected to any boundary cells as
\begin{equation}
  \label{eq:diff-iic}
  [\Delta_h \zeta_h]_i = [\Delta^2_h \psi_h]_i,\qquad i\in \IC\backslash \mathcal{BC}1.
\end{equation}
To obtain the discrete velocity field, one first maps the
stream function $\psi_h$ to a discrete scalar field $\tilde\psi_h$
defined at the cell vertices via \eqref{eq:psi-tilde}. The discrete
vector field representing the velocity is given by \eqref{eq:uh-def}. 

Substituting the above expressions for $q_i$, $\zeta_i$, $[\Delta_h
\zeta_h]_i$ and $\ub_h$ into the equation \eqref{eq:1a}, one obtains a
semi-discrete scheme,
\begin{multline}
\label{eq:semi-disc-iic}
  \dfrac{d}{d t}\left([\Delta_h \psi_h]_i + \beta y_i -
    \dfrac{f_0}{H}(\psi_i - b_i)\right)  + \left[ \nabla_h \cdot(\ub_h
    \hat q_h)\right]_i =\\ 
  \dfrac{1}{H}[\nabla\times\taub]_i - \alpha [\Delta_h \psi_h]_i +
  \mu[\Delta_h^2\psi_h]_i,\qquad 
  i\in\IC\backslash\BC1.
\end{multline}

We note that, in order to be able to explicitly enforce the boundary
conditions \eqref{eq:bc1a} and \eqref{eq:bc2a}, we need to treat the
(discrete) fourth order terms as unknowns. Thus, we discretize
\eqref{eq:semi-disc-iic}  by the 
semi-implicit Euler scheme, 
\begin{multline*}
  \dfrac{1}{\Delta t}\left([\Delta_h \psi^{n+1}_h]_i  -
    \dfrac{f_0}{H}\psi^{n+1}_i - [\Delta_h \psi^{n}_h]_i  +
    \dfrac{f_0}{H}\psi^{n}_i\right)   + \left[ \nabla_h \cdot(\ub^n_h
    \hat q^n_h)\right]_i =\\ 
  \dfrac{1}{H}[\nabla\times\taub^{n+1}]_i - \alpha [\Delta_h \psi^{n+1}_h]_i +
  \mu[\Delta_h^2\psi^{n+1}_h]_i\qquad 
  i\in\IC\backslash\BC1.
\end{multline*}
Moving all the unknowns (i.e.,~the $(n+1)$th-step unknown variables) to
the left-hand side of the equation, and the rest to the right-hand
side, we obtain
\begin{multline}
  \label{eq:imeu}
-\mu\Delta t\dfrac{g}{f_0}[\Delta_h^2\psi^{n+1}_h]_i + (1+\alpha\Delta
t) \dfrac{g}{f_0} [\Delta_h \psi^{n+1}_h]_i  - 
    \dfrac{f_0}{H}\psi^{n+1}_i  = \dfrac{\Delta
      t}{H}[\nabla\times\taub^{n+1}]_i -\\
    \Delta t\left[ \nabla_h \cdot(\ub^n_h
    \hat q^n_h)\right]_i + \dfrac{g}{f_0} [\Delta_h \psi^{n}_h]_i  -
    \dfrac{f_0}{H}\psi^{n}_i.
  \qquad 
  i\in\IC\backslash\BC1.
\end{multline}

We emphasize that the discrete equations \eqref{eq:imeu} are only for
interior cells not connected to any boundary cells. Additional
equations or constraints are needed to complete the system, and they
are supplied by discretizations of the boundary conditions
\eqref{eq:bc1a}, \eqref{eq:bc2a}, and the constraint
$\eqref{eq:bi-sys}_4$. The boundary condition \eqref{eq:bc1a} is
discretized as 
\begin{equation}
  \label{eq:imeu-bc1}
\psi^{n+1}_i = l,\qquad i\in\BC,
\end{equation}
where, for each time step, $l$ is an unknown scalar value. The BC
\eqref{eq:bc2a} can be discretized by specifying that the
stream function on the cells next to the boundary cells by assuming the same
scalar value as on the boundary cells,
\begin{equation}
  \label{eq:imeu-bc2}
  \psi^{n+1}_i = l,\qquad i\in\BC1.
\end{equation}
Finally, the constraint $\eqref{eq:bi-sys}_4$ is implemented by
requiring the area weighted sum of the stream function be zero,
\begin{equation}
  \label{eq:imeu-bc3}
  \sum_{i\in IC\cup \BC} A_i \psi_i = 0.
\end{equation}

To summarize, the semi-implicit numerical scheme for the viscous QG
equation consists of \eqref{eq:imeu}, \eqref{eq:imeu-bc1},
\eqref{eq:imeu-bc2}, and \eqref{eq:imeu-bc3}. The main advantage of
this approach is that the BC's are enforced explicitly and in a clean
way. The main disadvantage is that the discretization of the
fourth-order biharmonic operator results in a scheme with a stencil
that is 2-3 times larger than the stencil of a second-order elliptic
operator, making the scheme more expensive to solve.

\subsubsection{An explicit FV scheme for the viscous QG
  equation (VSFV2) }\label{s:fv-4} 
While the atmospheric and oceanic flows are considered viscous, the
viscosity is extremely small. For example, in the interior of the
world ocean, the flow has a Reynolds number on the order of $10^{20}$
(\cite{Cushman-Roisin2011-en}). Even though it is small, it would be a
mistake to simply ignore the viscosity, for numerical and dynamical
reasons. In fact, in the boundary layer region, both the length scale of
the flow and the Reynolds are much smaller than they are in the
interior of the domain, rendering the viscous terms dynamically
significant. On the other hand, as pointed out in the previous
subsection, incorporating the viscous terms implicitly in the solution
process can be costly. In this subsection, we introduce a new scheme
that is partly viscous and partly non-viscous. The viscosity is
included, {\itshape explicitly}, in the PV advection. But it is not
included when the stream function is computed from the PV. 

As before, we approximate the PV $q$ and the stream function $\psi$
with the cell-centered discrete scalar fields $q_h$ and $\psi_h$,
respectively. Once the PV $q_h$ is known, the stream function $\psi_h$
is computed from a discrete elliptic BVP,
\begin{equation}\label{eq:psi-q}
  \left\{
    \begin{aligned}
      &\dfrac{g}{f_0} [\Delta_h \psi_h]_i - \dfrac{f_0}{H} \psi_i = q_i
      - \beta y_i - \dfrac{f_0}{H}b_i, & & i\in \IC,\\
      &\psi_i = l & & i\in \BC,\\
      &\sum_{i\in \IC\cup\BC} \psi_i A_i = 0. & &
    \end{aligned}\right.
\end{equation}
We note that, while the continuous viscous QG equation requires two
boundary conditions, only one of the BC's,
i.e.~\eqref{eq:bc1a}, is explicitly enforced in the computation of
$\psi_h$. In other words, the flow is simply treated as inviscid as
far as the stream function is concerned. The stream function $\psi_h$ is
then mapped to a scalar field 
defined at cell vertices
\begin{equation}
  \label{eq:psi-tilde1}
  \tilde\psi_h = \tilde{[\psi_h]}.
\end{equation}
A discrete vector field representing the normal velocity components at
cell edges is computed from $\tilde\psi_h$ by
\begin{equation}
  \label{eq:uh1}
\ub_h = \tilde\nabla^\perp_h \tilde\psi_h.
\end{equation}

A discrete scalar field representing the vorticity is needed for
computing both the drag term and the diffusion term in
\eqref{eq:1a}. The discrete vorticity is computed as 
\begin{subequations}
\label{eq:zetah}
  \begin{align}
    \zeta_h &= \sum_{i\in\IC\cup\BC} \zeta_i
    \chi_i,\label{eq:zetah1a}\\
    \zeta_i &= [\Delta_h\psi_h]_i,\qquad i\in\IC\cup\BC. \label{eq:zetah1b}
  \end{align}
\end{subequations}
What distinguishes \eqref{eq:zetah} from \eqref{eq:zeta-ic} is the
range for the indices. The equation \eqref{eq:zeta-ic} is for interior
cells that are not connected to any boundary cells, while the formula
\eqref{eq:zetah} is for all cells. The specifications \eqref{eq:d34},
\eqref{eq:d26}, and \eqref{eq:d30} of relevant discrete differential
operators make it clear that the discrete Laplacian operator is
defined for all (interior and boundary) cells. It has been noted in
the comment following equation \eqref{eq:d26} that a no-flux condition
is implied on the boundary, because the edges along the domain
boundary are not included in the calculation of the discrete
divergence field. Hence, by computing the vorticity over the entire
domain with the formula \eqref{eq:zetah}, we implicitly enforce the
no-slip boundary condition \eqref{eq:bc2a}, 
\begin{equation*}
  \dfrac{\p\psi}{\p n} = 0,\qquad {\rm on}\quad\p\M,
\end{equation*}
which is required for the continuous viscous QG equation.

With the no-slip boundary condition implicitly enforced in the
calculation of the vorticity, it makes sense to update the PV on the
boundary
\begin{equation}
  \label{eq:pv-up}
  q_i = \zeta_i + \beta y_i - \dfrac{f_0}{H}\left(\psi_i - b_i\right),\qquad
  i\in \BC.
\end{equation}
The updated PV field is then used for computing the PV on the cell
edges
\begin{equation}
  \label{eq:qhat1}
  \hat q_h = \hat{[q_h]}.
\end{equation}

Finally, the PV evolves according to the semi-discrete equation 
\begin{equation}
\label{eq:q-visc}
  \dfrac{d q_i}{d t} + \left[ \nabla_h \cdot(\ub_h \hat q_h)\right]_i =
  \dfrac{1}{H}[\nabla\times\taub]_i - \alpha \zeta_i + \mu[\Delta_h
  \zeta_h]_i,\quad   i\in\IC,
\end{equation}
We note that the diffusion is included in the above explicitly. Since
the diffusion is only computed for the interior cells, no artificial
or non-physical boundary conditions are explicitly or implicitly
enforced on the vorticity.

\section{Analysis on numerical transport of the PV }\label{sec:anal}
We here deal
with the properties of an essential component of three out of the four
schemes presented in the previous section, namely the numerical
transport of the PV.  Rigorous numerical analysis of these schemes, concerning their
stability and convergence, will be left for future work. 

The PV plays a central role in geophysical fluid
dynamics (\cite{Lorenz2006-zx,Pedlosky1987-gk}). Its importance stems
from certain distinguishing properties that it enjoys. First,
according to \eqref{eq:1a}, in the absence of the external forcing and
the lateral and vertical diffusions, the PV is simply being
transported by the velocity field. In other words, in the absence of
the external forcing and the diffusions, the PV is preserved along the
fluid paths. This property is of vital importance in the
well-posedness analysis of the QG equation as well. 
Second, for predominantly two-dimensional large-scale geophysical
flows, it is well known that the energy tends to concentrate at large
scales (inverse energy cascade
(\cite{Batchelor1969-fo,Kraichnan1967-un,Lilly1969-wl}), and it is the
(potential) 
enstrophy, i.e.,~the second moment of the (potential) vorticity, that
cascades in the inertial range of the spectrum. 

In the following we provide an analysis of the performance of the PV advection
scheme on the PV dynamics and on the conservation of the potential
enstrophy. For our analysis, all the external influences, such as the
external forcing and the diffusions, have been dropped. 
Without these external factors, one will
see that our pure advection scheme preserves the PV along the fluid
paths in an averaging sense, and it also conserves the potential enstrophy
up to the time truncation errors. 

To evaluate the performance of the advection scheme on the PV
dynamics, we start with the scheme \eqref{eq:q-evol}. Substituting the
specification \eqref{eq:d26} for the discrete divergence operator in
the scheme and setting the external forcing $\taub$ and the bottom drag
$\alpha$ to zero, we obtain
\begin{equation}
  \label{eq:trspt1}
  \dfrac{d}{dt} q_i + \dfrac{1}{|A_i|}\sum_{e\in EC(i)} \hat q_e u_e
  l_e n_{e,i} = 0,\qquad i\in\IC\cup\BC.
\end{equation}
Using the identity $\hat q_e = \hat q_e - q_i + q_i$ in the above, we
have
\begin{equation}
  \label{eq:trspt2}
\dfrac{d}{dt} q_i + \dfrac{1}{|A_i|}\sum_{e\in EC(i)} (\hat q_e - q_i)
u_e l_e n_{e,i} + \dfrac{q_i}{|A_i|}\sum_{e\in EC(i)} u_e l_e n_{e,i}
= 0,\quad i\in\IC\cup\BC.\\
\end{equation}
We note that the last term on the left-hand side of \eqref{eq:trspt2}
vanishes because it involves the discrete divergence of $\ub_h$, which
is divergence-free on 
every cell (\cite{Chen2016-gl}). We also note that $\hat q_e$,
computed via \eqref{eq:qhat_e}, is the arithmetic mean of $q_i$ and the
value of the PV on the cell on the other side of edge $e$, and $\hat
q_e - q_i$ will simply be half of the difference between these two
discrete PV values, and oriented outwards from cell $i$, i.e.
\begin{equation}
  \label{eq:trspt3}
  \hat q_e - q_i = \dfrac{1}{2} d_e [\nabla_h q_h]_e n_{e,i}.
\end{equation}
Substituting \eqref{eq:trspt3} in \eqref{eq:trspt2}, we obtain
\begin{equation}
  \label{eq:4}
  \dfrac{d}{dt} q_i + \dfrac{1}{|A_i|}\sum_{e\in EC(i)}
  \dfrac{1}{2}d_e l_e [\nabla_h q_h]_e u_e = 0,\qquad i\in\IC\cup\BC.
\end{equation}
Noticing that, for the diamond-shaped region,
\begin{equation*}
  A_e = \dfrac{1}{2} d_e l_e,
\end{equation*}
we can write \eqref{eq:4} as
\begin{equation}
  \label{eq:trspt5}
\dfrac{d}{dt} q_i + \sum_{e\in EC(i)} \dfrac{|A_e|}{|A_i|}[\nabla_h
q_h]_e u_e = 0,\quad i\in\IC\cup\BC.
\end{equation}
Within the summation, $[\nabla_h q_h]_e u_e$ is obviously an
approximation of the inner product of the {\itshape normal} components
of $\nabla q$ and $\ub$ on the edge. For the centroidal Voronoi tessellation meshes (\cite{Du2003-gn,Du1999-th}) , two
neighboring cell centers are equi-distant to the edge between them,
and, thus, the following identity holds,
\begin{equation*}
  |A_i | = \dfrac{1}{2} \sum_{e\in\EC(i)} |A_e|\qquad\textrm{or}\qquad
  \sum_{e\in\EC(i)} \dfrac{|A_e|}{|A_i|} = 2. 
\end{equation*}
Therefore, the summation on the left-hand side of \eqref{eq:trspt5} is
a linear convex combination of the discrete convections  across the
edges surrounding cell $i$. The factor 2 comes from the fact that only
the {\itshape normal} components are represented in the linear
combination. 

In the absence of an external forcing and the diffusions, the
potential enstrophy is conserved under the QG equation. This can be
directly verified by multiplying \eqref{eq:1a} with $q$ and integrate
over $\M$, 
\begin{align*}
  \dfrac{d}{dt}\int_\M q^2 dx 
&= -\int_\M \ub\cdot\nabla\dfrac{1}{2}q^2 dx \\
&= -\int_\M \nabla\cdot(\dfrac{1}{2}\ub q^2)dx \\
&= -\int_{\p\M} \dfrac{1}{2}q^2 \ub\cdot\nb ds \\
&= 0. 
\end{align*}

We now show that our PV convection scheme also conserves the potential
enstrophy, up to the truncation errors in time. We again start with
the scheme \eqref{eq:q-evol}, and drop the external forcing and the
bottom drag,
\begin{equation}
  \label{eq:trspt6}
  \dfrac{d}{d t} q_h + \nabla_h \cdot(\ub_h\hat q_h) = 0.
\end{equation}
Taking the inner product of \eqref{eq:trspt6} with $q_h$, we have
\textcolor{rev}{\begin{equation*}
  \dfrac{d}{d t}(q_h, q_h)_{0,h}  + (\nabla_h \cdot(\ub_h\hat q_h), \, q_h)_{0,h} = 0.
\end{equation*}
In the above, both $(\cdot,\cdot)_{0,h}$ designate the inner product
\eqref{eq:inp1} 
of scalar fields on the primary mesh.}  Applying the discrete
integration by parts formula \eqref{eq:d41} to 
the second term, we obtain
\textcolor{rev}{\begin{equation}
  \label{eq:trspt7}
  \dfrac{d}{d t}(q_h, q_h)_{0,h}  - 2(\ub_h\hat q_h, \, \nabla_h q_h)_{0,h} = 0.
\end{equation}
Here, the second $(\cdot,\cdot)_{0,h}$ designates the inner product
\eqref{eq:inp2} of
discrete normal vector fields on the edges.}
We note that $\hat q_h = \sum \hat q_e \chi_e$, and $\hat q_e$ is the
arithmetic mean of the discrete PV values on the neighboring cells on
edge $e$, and $\nabla_h q_h = \sum [\nabla_h q_h]_e \chi_e\nb_e$,
where $[\nabla_h q_h]_e$ is the gradient between the discrete $PV$
values on the neighboring cells on edge $e$. It is easy to see that,
\begin{equation*}
  \hat q_h \nabla_h q_h = \sum_{e\in\IE\cup\BE} \dfrac{1}{2}[\nabla_h
  q^2_h]_e  \chi_e \nb_e\equiv \nabla_h \dfrac{1}{2} q^2_h.
\end{equation*}
Therefore, we can write \eqref{eq:trspt7} as 
\begin{equation}
  \label{eq:trspt8}
  \dfrac{d}{d t}(q_h, q_h)_{0,h}  - (\ub_h, \, \nabla_h q^2_h)_{0,h} = 0.
\end{equation}
Applying the integration by parts formula \eqref{eq:d41} again, we
arrive at 
\begin{equation*}
  \dfrac{d}{d t}(q_h, q_h)_{0,h}  + (\nabla_h\cdot \ub_h, \,  q^2_h)_{0,h} = 0.
\end{equation*}
Thanks to the incompressibility of the discrete vector field $\ub_h$, we
have
\begin{equation*}
  \dfrac{d}{d t}(q_h, q_h)_{0,h}  = 0.
\end{equation*}
The potential enstrophy is conserved up to the time discretization errors.

\section{Numerical experiments}\label{s:numerics}

\subsection{A test case with a freely-evolving circular
  flow}\label{sec:circ}

In the previous section we have shown that the transport part of our
scheme {\itshape conserves the PV along the fluid path { in the average
  sense}, and conserves the potential enstrophy up to the time
truncation error}. In order to evaluate the performance of the scheme in
practical simulations where space averaging and time discretization
errors are present, we apply the scheme to a freely-evolving
circular flow with no external forcing, no diffusion or no bottom
drag. 

\begin{figure}[!ht]
  \centering
  \includegraphics[width=3.5in]{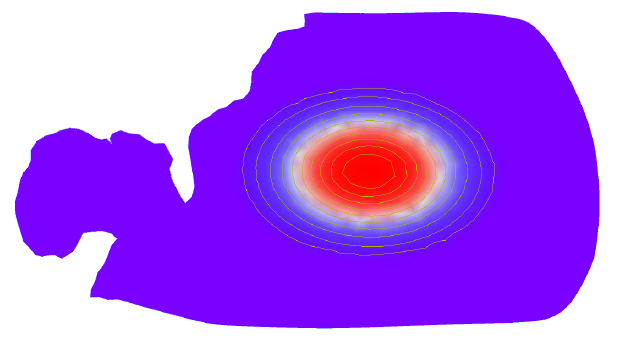}
  \caption{The initial stream function for a test of the pure PV
    advection scheme.}
  \label{fig:trspt-test}
\end{figure}

We consider one section of the mid-latitude northern Atlantic ocean. 
The initial state of the flow is given by 
\begin{align*}
   d &= \sqrt{\dfrac{(\theta -
      \theta_c)^2}{\Delta\theta^2} + \dfrac{(\lambda -
      \lambda_c)^2}{\Delta\lambda^2}} \\
  \psi(\lambda,\theta) &= e^{-d^2} \times (1- \tanh(20*(d-1.5))
\end{align*}
In the above, $\theta_c = 0.5088$ and $\lambda_c = -1.1$ are the
latitude and longitude of the center point of the Atlantic section,
$\Delta\theta = .08688$, $\Delta\lambda = .15$. The $\tanh$ function
is used to ensure that the initial stream function is flat (with one constant
value) along the boundary of the domain. A plot of the domain and of
the stream function is shown in Figure \ref{fig:trspt-test}. This
stream function produces a circular clockwise velocity field in the
middle of the ocean, with the
maximum speed of about $.8$m/s near the center, which is considered
fast for oceanic flows. No external forcing or diffusion is applied,
and the flow is allowed to evolve freely. The circular pattern will
breakdown eventually, because of the irregular shape of the domain and
because of the non-uniform Coriolis force.

A spherical centroidal Voronoi-Delaunay mesh of the domain with 186,168 primary Voronoi cells and
367,589 dual Delaunay triangles is used, and the mesh resolution gradually  changes from 21km to 2km when 
getting closer and closer to the boundary.
 \textcolor{rev}{The numerical scheme (IVFV1) of Section \ref{s:fv-1}
   is used, but with the external 
forcing and the bottom drag $\alpha$ set to zero so that only the PV
transport participates in the dynamics.} A fourth-order Runge-Kutta
scheme is used to discretize the equation \eqref{eq:q-evol}. 
The model is run for 10 years, and the  time step size is $\Delta t =
1350s$.  

\begin{figure}[!ht]
  \centering
  \includegraphics[width=7.3in]{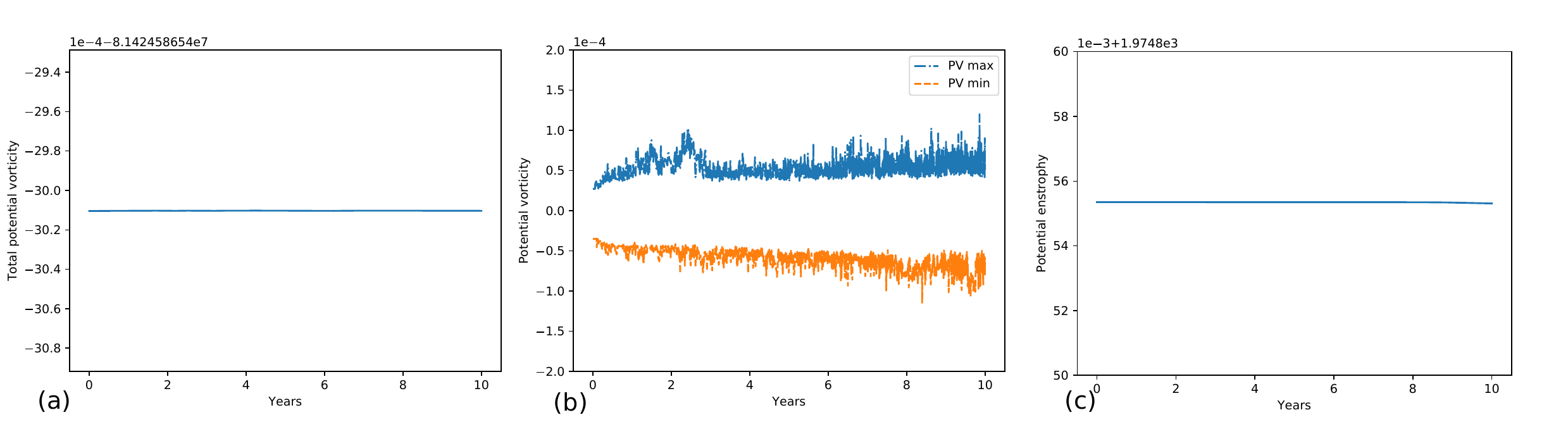}
  \caption{(a) Conservation of the total potential vorticity to the
    machine error. Over the course of 10 years, the change in the
    total potential vorticity is 
    about $1.25\times 10^{-16}$ of the original total PV. (b) The
    evolution of the maximum and minimum values of the PV. (c) The
    conservation of the potential enstrophy. Over the course of 10
    years, the change in the potential enstrophy is about $3\times
    10^{-8}$ of the original potential enstrophy. 
  }
  \label{fig:pv-conserv}
\end{figure}

For this simulation, we want to see how well the PV and the potential
enstrophy are conserved. For the QG equation, since the velocity field
is incompressible, the PV satisfies the conservation laws. Due to the
finite volume discretization technique used for the scheme, the PV is
expected to be conserved both locally and globally. The global
conservation of the PV is verified
by a plot of the total PV against time (Figure \ref{fig:pv-conserv}
(a)). The total PV remains constant in time. Over the course of 10
years, the change in the total PV is only about $1.25\times 10^{-16}$
of the original total PV. Thus, under this scheme, the total PV is
conserved up the machine errors. 

In the continuous system, the PV is also a quantity that is being
transported by the velocity field, which means that, in the absence of
an external forcing or diffusion, the PV is conserved along the fluid
path. In the previous section, it has been shown that the transport
component of our schemes conserve the PV along the fluid path in the
averaging sense. From the discrete velocity field, it is impossible to
fully reconstruct the fluid path, and therefore it is impossible to
observe how PV evolves along any given fluid paths. Instead, we focus
on the maximum value and the minimum values of the PV field. In the
continuous system, since the PV is just being transported along the
fluid particles, its maximum and minimum values will remain constant
for all time. We want to know how well our scheme is able to preserve
the maximum and minimum values of the PV field. However, it is not
reasonable to expect the maximum and minimum values of the discrete PV
are preserved exactly, because the scheme has only been shown to
conserve the PV along the fluid path in the averaging
sense. Oscillations are expected in the PV values, especially on the
extreme  values. 
In Figure
\ref{fig:pv-conserv}(b), we plot the maximum and minimum values of the
discrete PV against time. It can be seen the maximum and minimum values
of the PV are not preserved exactly. They oscillate very fast, but
they remain within some well-defined bounds. Initially, the PV values
are within the interval 
of $[-0.3\times 10^{-4}, 0.3\times 10^{-4}]$, and throughout the entire
simulation period of 10 years, the PV values remain within the bound
of $[-1.2\times 10^{-4}, 1.2\times 10^{-4}]$. 

Finally, we examine the performance of the scheme in conserving the
potential enstrophy, which is the second moment of the PV. In Figure
\ref{fig:pv-conserv}(c), we plot the 
potential enstrophy of the discrete system against time. The plot shows
that the potential enstrophy remains largely constant in time. An
close examination of the numerical values shows that, over the course
of 10 years, the change in the potential enstrophy is about $3\times
10^{-8}$ of the initial potential enstrophy, which confirms that,
under this scheme, the potential enstrophy is conserved up to the time
discretization errors.

\subsection{A test case with the intensified western boundary
  current}\label{s:gulf}
In this subsection, we numerically study the qualitative behaviors of
all the schemes presented in Section \ref{sec:scheme}, and identify
the scheme that is most suitable for simulation of realistic flows in
terms of both performance and efficiency.

To accomplish these goals,
we utilize a test case that has much relevance to the realistic
oceanic flows, namely a model for the intensified western boundary
current, which, in the real world, takes the form of the Gulf stream,
the Kuroshio current, the Brazilian current, etc.  We consider the
same section of the Atlantic ocean as in the previous section. But,
here, the flow is initially at rest, and is being constantly driven by
a zonally uniform wind stress of the form
\begin{equation*}
  \tau = \tau_0 \cos\left( \dfrac{\pi(\theta - \theta_0)}{\Delta
      \theta}\right), 
\end{equation*}
where $\tau_0 = 1.0\times 10^{-4}\textrm{ m}^2/\textrm{s}^2$, $\theta$
the latitude, $\theta_0 \approx 15^\circ$ the latitude of the southern
edge of the domain, and $\Delta \theta\approx 30^\circ$ the latitudinal
width of the domain. The wind stress has a negative curl, and thus is a
source of negative vorticity to the flow.

\begin{figure}[!ht]
  \centering
  \includegraphics[width=3.4in]{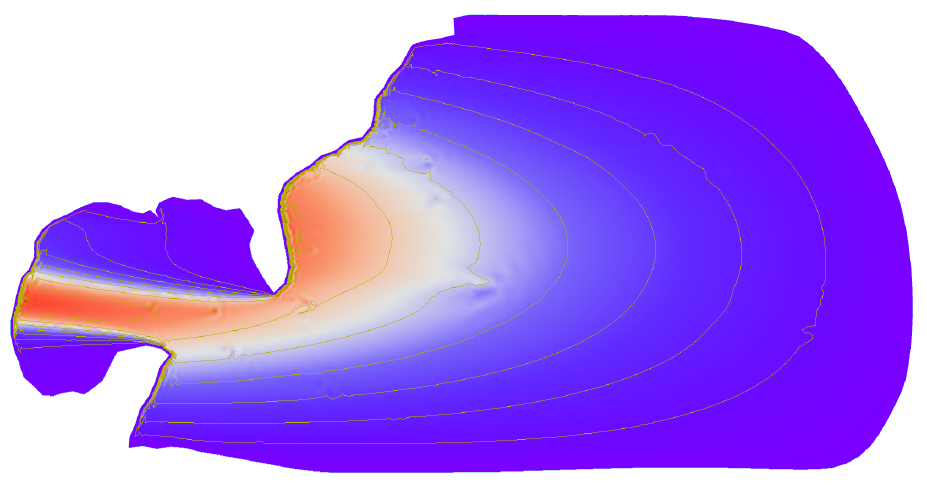}
  \caption{The solution (stream function)  to the steady state linearized equation with a
    bottom drag $\alpha = 5\times 10^{-8}$ (the Stommel boundary layer) . }
  \label{fig:stommel}
\end{figure}

\begin{figure}[!ht]
  \centering
  \includegraphics[width=3.4in]{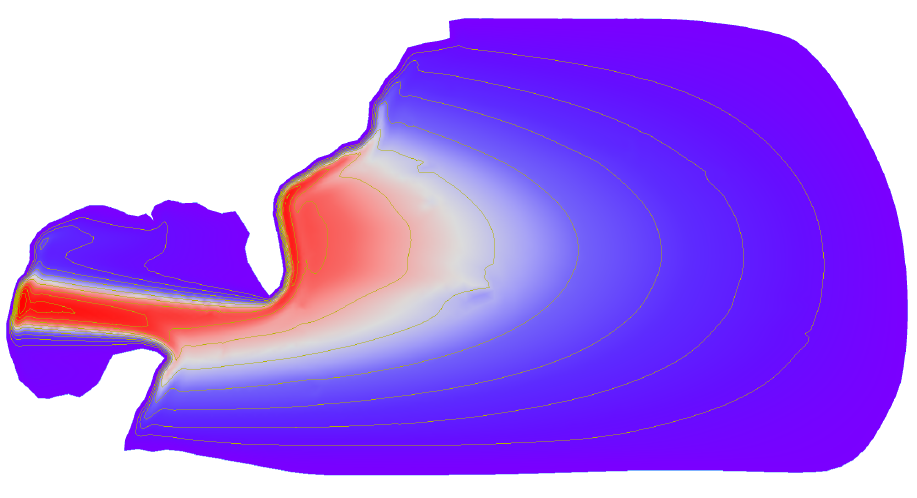}
  \caption{The solution (stream function) to the steady state linearized equation  with
    a viscosity $\mu = 100$ (the Munk boundary layer). }
  \label{fig:munk}
\end{figure}

Stommel proposed the very first model (\cite{Stommel1948-wx}) for the
Gulf stream, articulating that the Gulf Stream is primarily the
consequence of 
the anisotropic effect of the Coriolis force. What is the most
surprising about Stommel's model is that the bottom drag, taking the
same form as the zeroth order term on the right-hand side of
\eqref{eq:1a},  is the only 
damping term, even though the Gulf stream has the appearance of a
lateral boundary layer. When the lateral diffusion is included, the
resulting boundary layer is called the Munk layer. Both the Stommel
model and the Munk model are  steady-state and linear. The solutions
to these steady-state and linear models over the domain under
consideration  and with the wind stress are plotted in
Figures \ref{fig:stommel} and \ref{fig:munk}. These plots clearly show
that the Stommel layer and the Munk layer share the same structure;
they differ only in details near the western boundary. These
steady-state solutions will serve as the reference solutions for the
time-dependent solutions that we are about to compute. 

\begin{figure}[!ht]
  \centering
  \includegraphics[width=3.4in]{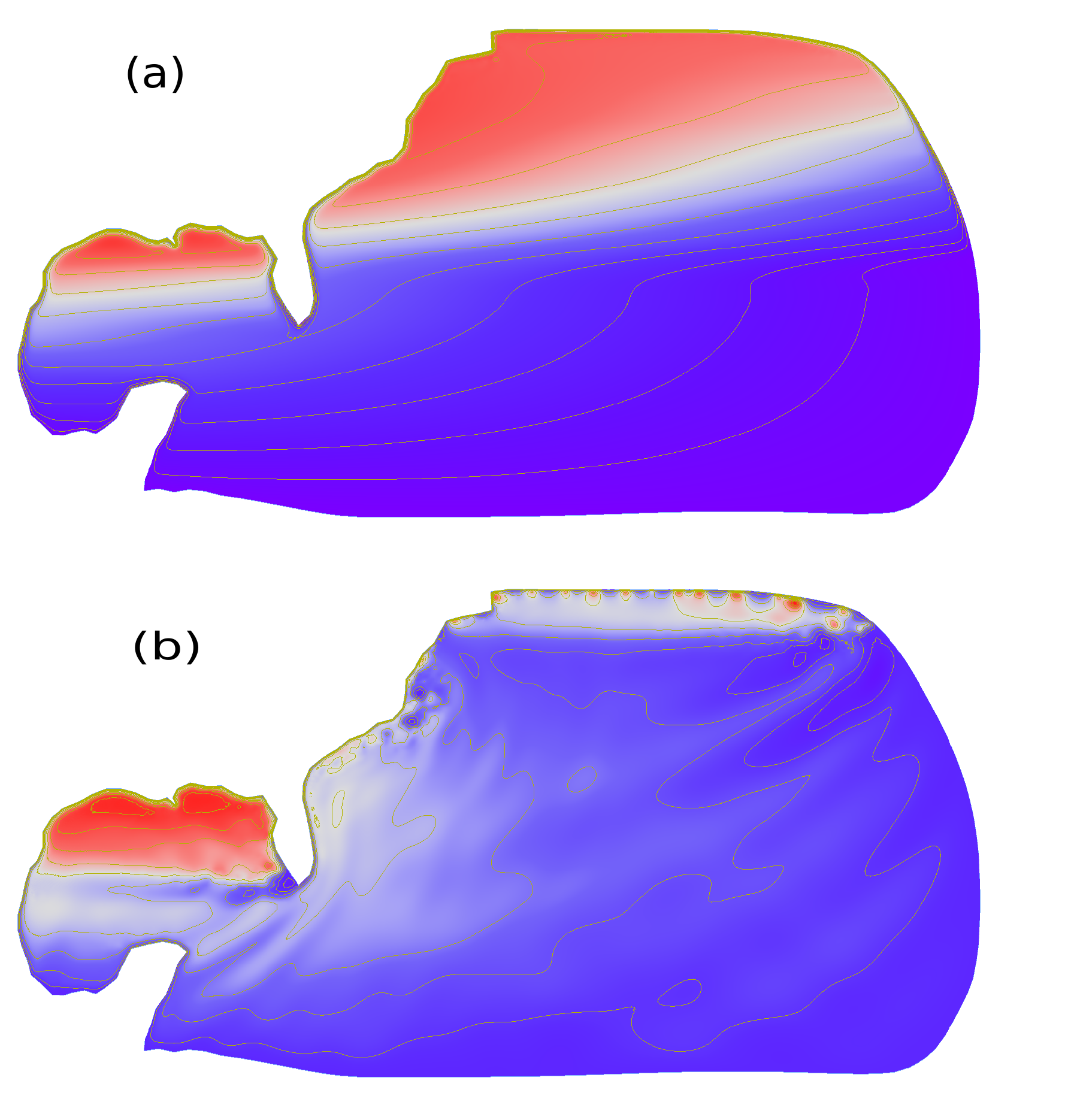}
  \caption{The stream functions at year 5, as computed with the
    FV schemes for the inviscid models. The flows are dissipated with
    only  a bottom drag of 
    $\alpha = 5\times 10^{-8}$. Top: IVFV1; bottom: IVFV2.}
  \label{fig:y5-invi}
\end{figure}

In the inviscid model, the bottom drag is the only dissipative
mechanism. In the previous Section, two schemes were proposed for this
model. The scheme IVFV1 of Section \ref{s:fv-1} only enforces the BC
\eqref{eq:bc1a} theoretically required for the model,  while the
scheme IVFV2 of Section \ref{s:fv-2} enforces the BC \eqref{eq:bc1a} and an
extra artificial BC \eqref{eq:psi-bc-2}. Both schemes are applied to
this test problem, with the bottom drag parameter $\alpha = 5\times
10^{-8}$. The stream function computed using the scheme IVFV1, at the
end of the 5-year simulation, is shown in Figure
\ref{fig:y5-invi} (a). The domain is unmistakably divided into two
regions, one with positive vorticity, and the other with flat
stream function and (nearly) zero vorticity. In addition to the western
boundary, the flow also exhibits intensified boundary currents at the
northern and eastern boundaries. While this result is not a good
reflection of the realistic oceanic flows, where the intensified
boundary currents are only observed on the western boundary, it does
confirm an important principle for large-scale geophysical flows,
namely that, in the absence of frictions, the PV tend to homogenize
its self in the interior of the flow (\cite{Rhines1982-te}). The
region with negative vorticity on the top is trying to offset the
increase in the planetary vorticity in the flow. 

With the artificial boundary condition \eqref{eq:psi-bc-2} imposed by
the scheme IVFV2, the homogenization process is disrupted. Since the
relative vorticity is set to zero on the boundary and the surface
fluctuations remain small, the PV on the boundary remain close to the
planetary vorticity. This is approximately the case in reality (the
planetary vorticity always dominates), and as a consequence, the
stream function (Figure \ref{fig:y5-invi} (b)) is also much closer to those
in Figures \ref{fig:stommel} and \ref{fig:munk}. The currents are only
intensified near the western boundary. The stream of vortices on the
top of the domain is the result of the combination of the free slip
boundary condition and the vanishing curl of the wind stress there.

\begin{figure}[!ht]
  \centering
  \includegraphics[width=3.4in]{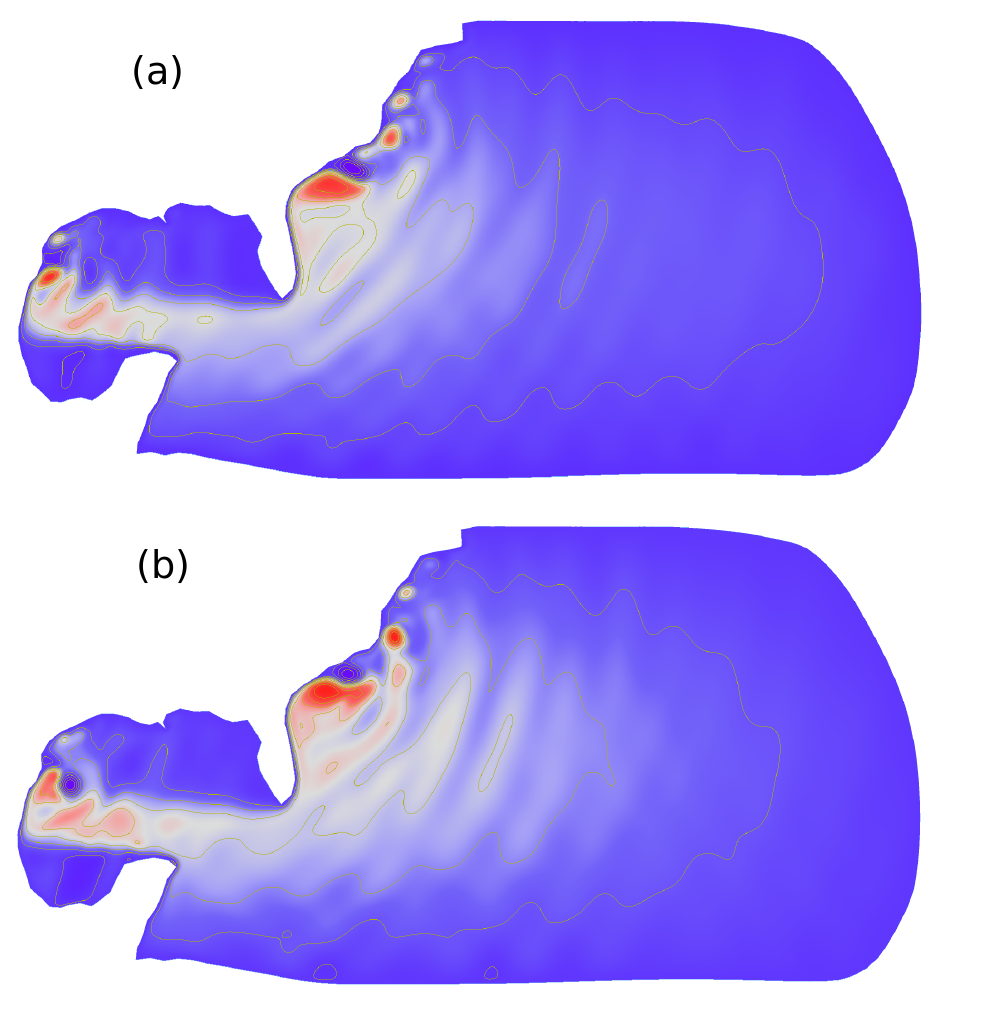}
  \caption{The stream functions at year 5, as computed with the FV
    schemes for the viscous models. The flows are dissipated with  a
    bottom drag of 
    $\alpha = 3\times 10^{-8}$  and  a viscosity $\mu = 40$. Top: VSFV1; bottom: VSFV2.}
  \label{fig:y5-visc}
\end{figure}

\begin{figure}[!ht]
  \centering
  \includegraphics[width=5in]{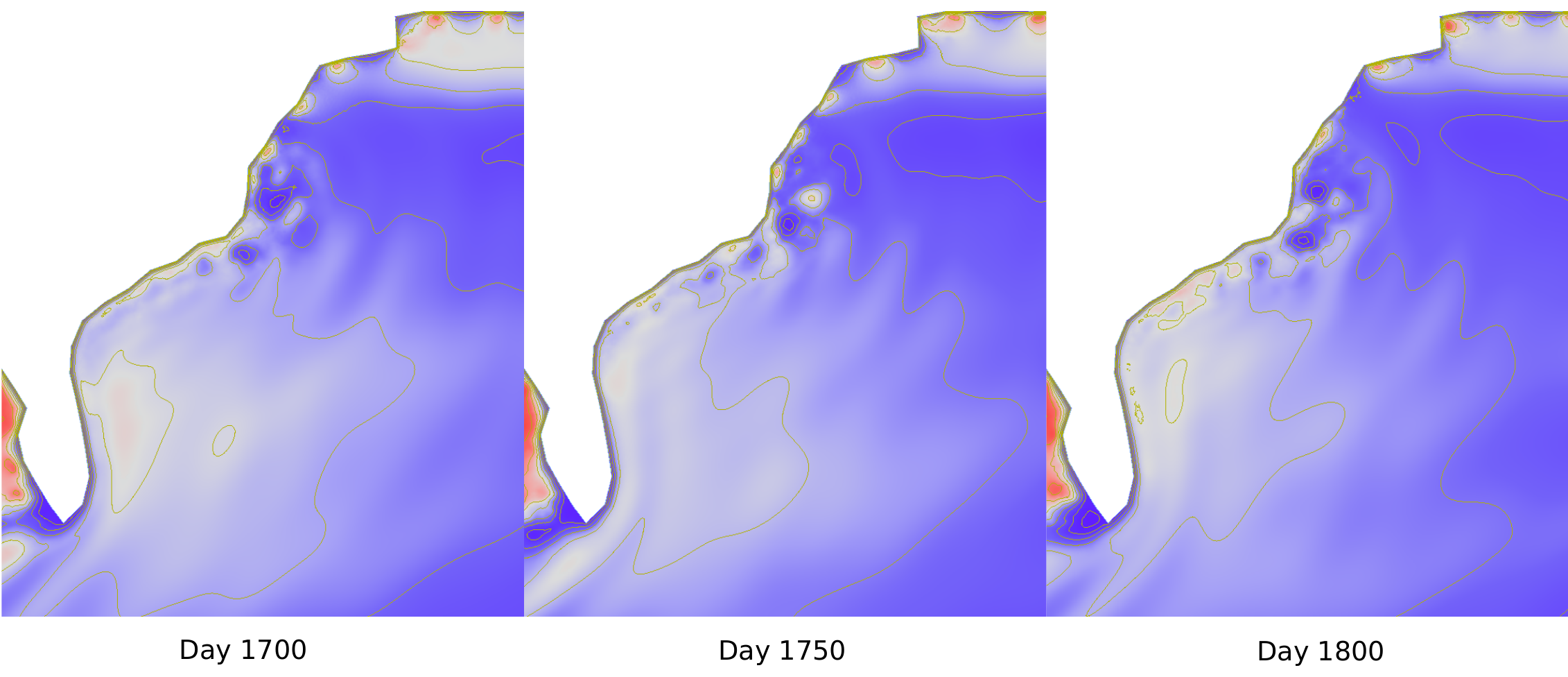}
  \caption{The stream functions, at days 1700, 1750, and 1800, with a bottom drag $\alpha
  = 3\times 10^{-8}$,
    computed using IVFV1  and the Runge-Kutta 4th 
    order scheme. }
  \label{fig:psi-zoom-bottom}
\end{figure}

\begin{figure}[!ht]
  \centering
  \includegraphics[width=5in]{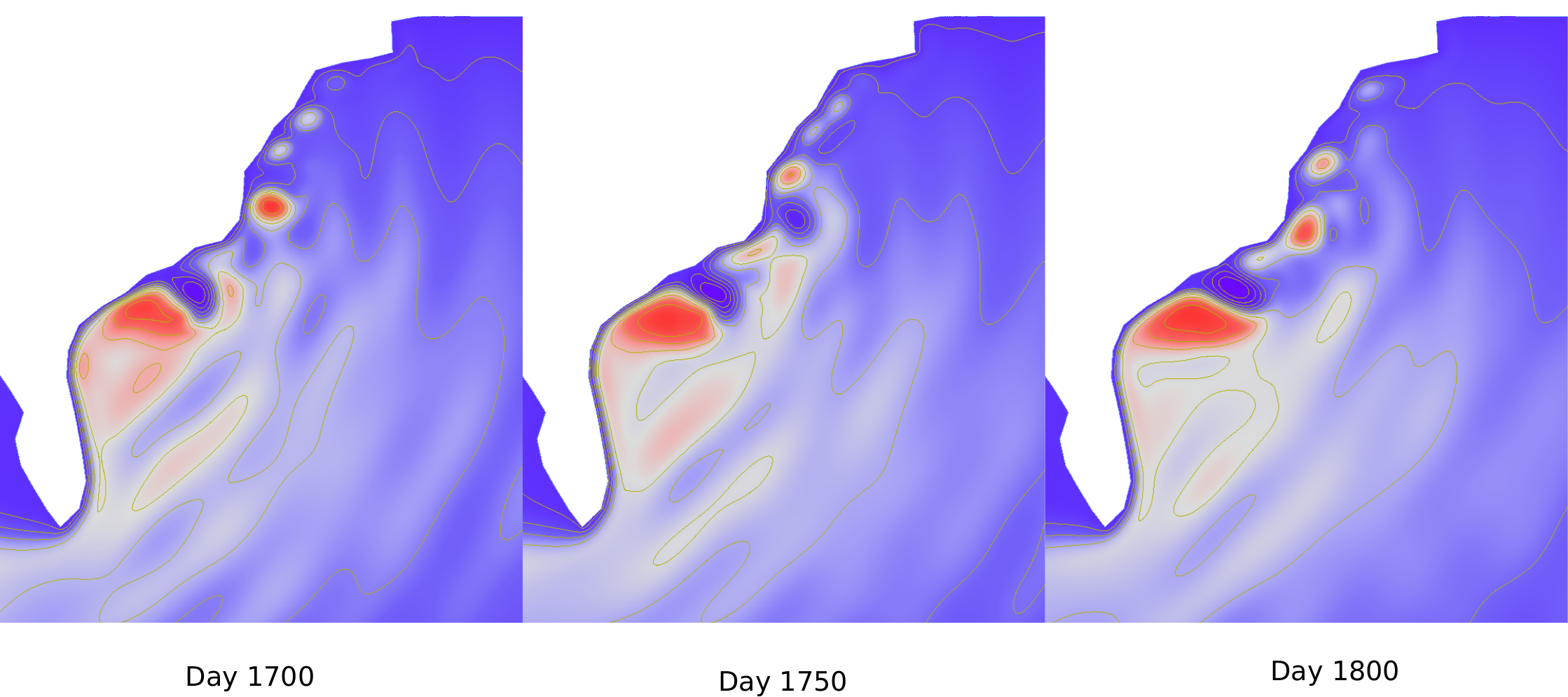}
  \caption{The stream functions, at days 1700, 1750, and 1800, with a bottom drag $\alpha
  = 3\times 10^{-8}$  and a viscosity $\mu = 40$,
    computed using VSFV1 and the semi-implicit Euler scheme. }
  \label{fig:psi-zoom-be}
\end{figure}

\begin{figure}[!ht]
  \centering
  \includegraphics[width=5in]{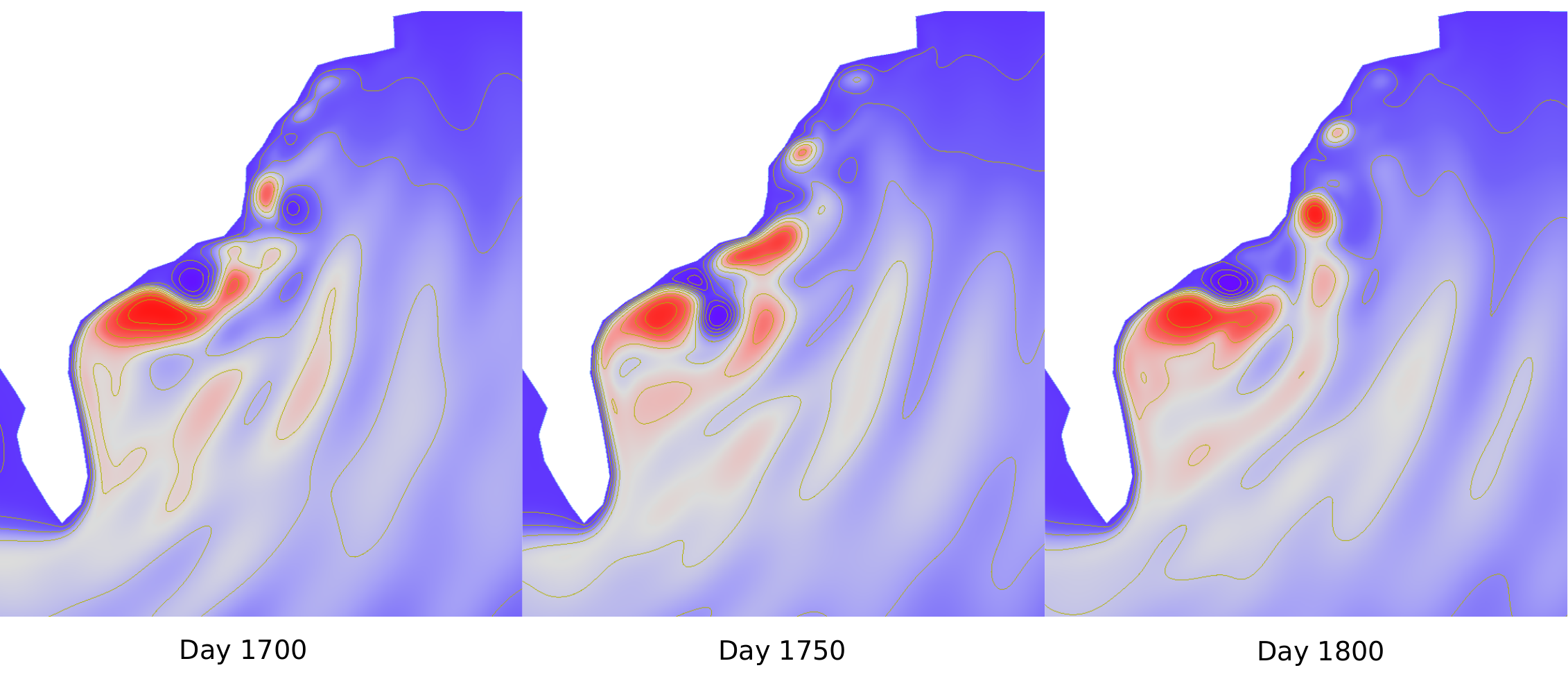}
  \caption{The stream functions, at days 1700, 1750, and 1800, with a bottom drag $\alpha
  = 3\times 10^{-8}$  and a viscosity $\mu = 40$,
    computed using VSFV2  and the Runge-Kutta 4th 
    order scheme. }
  \label{fig:psi-zoom-rk}
\end{figure}

When the lateral diffusion is included in the viscous model, the extra no-slip BC
\eqref{eq:bc2a} is required for the continuous model. Both the BC's
\eqref{eq:bc1a} and \eqref{eq:bc2a} are enforced explicitly in the
semi-implicit scheme VSFV1, whereas the BC \eqref{eq:bc2a} is only
implicitly enforced in the explicit scheme VSFV2. Both schemes are
applied to this case, with a bottom drag $\alpha = 3\times 10^{-8}$,
and the lateral diffusion $\mu = 40$. The stream functions at the end
of the 5-year simulation are plotted in Figure \ref{fig:y5-visc} (a)
and (b). The stream functions from both schemes share the same
large-scale structure of the steady-state Stommel-Munk boundary
layers. Indeed, both schemes produce an intensified current near the western
boundary only. The stream of vortices also disappear from the top of
the domain, thanks to the no-slip boundary condition. 

A comparison of the details in the boundary layer region is warranted
in order to further distinguish these schemes. 
Presented in Figures \ref{fig:psi-zoom-bottom}, \ref{fig:psi-zoom-be}
and \ref{fig:psi-zoom-rk} are detailed views of the boundary layer
region along the east coast of the United States for the schemes
IVFV2, VSFV1, and VSFV2, respectively. IVFV1 is not included in this
comparison because it does not produce a realistic rendering of the
intensified western boundary current. These figures reveal several
aspects about these schemes. First, while the artificial BC
\eqref{eq:psi-bc-2} helps the scheme IVFV2 produce a stream function
field that is close to that of the Stommel-Munk boundary layer, it
also makes the boundary layer thinner than those produced by the
scheme VSFV1 and VSFV2. Second, the artificial BC \eqref{eq:psi-bc-2}
also leads to subdued  eddy activities in the boundary layer region,
compared with those in Figure \ref{fig:psi-zoom-be} and
\ref{fig:psi-zoom-rk}. Lastly, while the schemes VSFV1 and VSFV2 impose
the theoretical BC \eqref{eq:bc2a} in different fashions, the boundary
layers that they produce are of similar structures, with comparable
width and eddy activity levels.

While the schemes VSFV1 and VSFV2 have similar dynamical behaviors,
the scheme VSFV2, coupled with a higher-order explicit scheme such the
the 4th order Runge-Kutta method, has significant efficiency advantages over the scheme
VSFV1. The first advantage comes from the fact that the scheme VSFV2
only has to solve a second-order elliptic 
problem, while VSFV1 has to solve a fourth-order elliptic problem,
which has a larger discrete stencil. 
Our simulation results  indicate that solving a
fourth-order elliptic PDE is more than twice as expensive as solving a
second-order elliptic PDE. Another surprising and much more
significant advantage of the scheme VSFV2 is that, with the 4th
order Runge-Kutta 
method, it can use a comparable or larger time step size than the
scheme VSFV1. This advantage can be attributed to the fact that the
fourth-order Runge-Kutta method has a very large stability region. In our numerical tests, the scheme VSFV2 runs smoothly with a time step size
of even 10800s, while the scheme VSFV1 can only use a time step size of 1200s. 
While using a different higher-order implicit time stepping technique
may be able to ameliorate the efficiency for VSFV1, it is not expected
that it 
can achieve the same performance as VSFV2. 

\section{Conclusions}\label{s:conclu}

The current work is part of the project to develop
vorticity-divergence based numerical schemes for large-scale
geophysical flows on arbitrarily unstructured primal-dual meshes over bounded
domains. Vorticity-divergence based numerical schemes offer the
optimal representations of the inertial-gravity wave relations, which
are essential for an accurate simulation of the geostrophic adjustment
process. They advance purely scalar
quantities, i.e.~the mass, vorticity, and divergence, and thus avoid
the need for vector reconstructions, for which there is no
satisfactory solutions yet on unstructured meshes. Finally,
vorticity-divergence based numerical schemes, through the
finite-volume finite-difference discretization techniques,  offer the
opportunity for a
tighter control over such fundamental quantities as vorticity, and
thus a more accurate simulation of the overall dynamics. 

The current work concerns the barotropic QG equation with a
free surface, which is a relatively simple model on the hierarchy of
models for geophysical flows. Several finite volume schemes are
presented for both the  inviscid
and the viscous equations. Great care is taken with regard to the
specification of the boundary conditions. For the inviscid equation, a
no-flux BC is required for the well-posedness of the model
(\cite{Chen2017-fh}). 
We present a scheme (IVFV1) that implements this BC explicitly. But
this scheme fails to produce a realistic rendering of the intensified
western boundary current, due to the lack of control over the PV on
the boundary. A second scheme (IVFV2) with an artificial free-slip BC
is also presented. This scheme adds a certain level of control on the
PV on the boundary by mandating that the vorticity be zero there. The
rendering of the intensified western boundary current is improved,
but, due to the free-slip condition, an unphysical stream of vortices
appear along the northern boundary. The viscous model requires both
no-flux and no-slip boundary conditions. Two schemes are proposed for
this model as well. The scheme VSFV1 implements both the no-flux and
no-slip BCs explicitly by using an semi-implicit time stepping
scheme. The VSFV2 implements the no-flux BC explicitly, but the
no-slip BC implicitly. Both schemes produce qualitatively similar
results in the test involving the intensified western boundary
currents. We find that the scheme VSFV2, coupled with the 4th order
Runge-Kutta scheme, is significantly more efficient than the VSFV1
scheme. Realistic geophysical flows typically involve the lateral
diffusion, usually in the form of eddy viscosity. The tests performed
in this study show that, in these cases, 
the scheme VSFV2 should be used, as it can  deliver the
right dynamics at relatively little cost (compared with VSFV1). 
The scheme IVFV1, thanks to its clean implementation of the BC, is of
theoretical interest, as it can be used in studies focusing on the
inviscid nonlinear dynamics of the geophysical flows. 

In subsequent works, we will consider more complex systems, such as
the shallow water equations, the isopycnal models, etc. These system
involve more prognostic and diagnostic variables, and it will be more
challenging to  find the right
specifications for these variables so that certain quantities, e.g.~the
PV, the potential enstrophy etc., are conserved. 
As mentioned in the introduction, in order for the vorticity-divergence based
schemes to be appealing for practical applications, the performance of
these schemes need to be further improved. The performance issue was skipped
over all together in the current work, but will be dealt with in
subsequent works. Our strategy for dealing with this issue will be a
combination of iterative solvers, preconditioners, and good initial
guesses. 

\begin{acknowledgements}
 This work is partially supported by the Simons Foundation
(contract number~319070) and by the US
Department of Energy (grant number~DE-SC0016540). We are also grateful
to the referees whose comments helped improve the paper. 
\end{acknowledgements}

\bibliographystyle{spmpsci}      
\bibliography{paperpile}   

\end{document}